\documentclass[preprint,1p]{elsarticle}
\usepackage{amsmath,latexsym,amssymb}
\usepackage[utf8]{inputenc}
\usepackage[T1] {fontenc}
\usepackage{lmodern}
\newif\ifpdf
\ifx\pdfoutput\undefined
\pdffalse 
\else
\pdfoutput=1 
\pdftrue
\fi

\ifpdf
\else
\fi

\newtheorem{lemma}{Lemma}[section]
\newtheorem{theorem}[lemma]{Theorem}
\newtheorem{corollary}[lemma]{Corollary}
\newtheorem{definition}[lemma]{Definition}

\newtheorem{remark}[lemma]{Remark}
\newtheorem{proposition}[lemma]{Proposition}

\newtheorem{conjecture}[lemma]{Conjecture}
\newtheorem{example}[lemma]{Example}

\newtheorem{question}[lemma]{Question}
\newenvironment{proof}[1][Proof]{\textbf{#1.} }{\ \rule{0.5em}{0.5em}}

\usepackage{setspace}

\begin{document}

\ifpdf
\DeclareGraphicsExtensions{.pdf, .jpg, .tif}
\else
\DeclareGraphicsExtensions{.eps, .jpg}
\fi

\title{Distance-transitive digraphs: descendant-homogeneity, property $Z$ and reachability\tnoteref{t1}}
\tnotetext[t1]{This work was supported by FAPDF.} 

\author{Daniela A. Amato\corref{cor}}
\ead{amatodani@gmail.com, d.a.amato@mat.unb.br}

\address{Departamento de Matem\'atica, Universidade de Bras\'ilia,  Campus Universit\'ario Darcy Ribeiro, CEP 70910-900, Brazil.}

\date{}

\begin{abstract} 
A digraph is {\itshape distance-transitive} if for every integer $s \geq 0$, the automorphism group of the digraph is transitive on pairs $(u,v)$ of vertices for which there is a directed path of length $s$ from $u$ to $v$, but no directed path of length $t<s$. This implies vertex and edge transitivity but is weaker than being highly-arc-transitive. A digraph is said to have {\itshape property Z} if it has a homomorphism onto the two-way infinite directed path. In a digraph $D$, an edge $e'$ is {\itshape reachable} from an edge $e$ if there exists an alternating walk in $D$ whose initial and terminal edges are $e$ and $e'$. Reachability is an equivalence relation on the edge set of $D$ and if $D$ is transitive on edges, then this relation is either universal or all of its equivalence classes induce isomorphic bipartite digraphs.

We investigate the class of  infinite distance-transitive digraphs of finite out-valency. First, we show that earlier results, proved in the context of highly-arc-transitive digraphs, hold for the class of distance-transitive digraphs. Second, we show that if $D$ is a weakly descendant-homogeneous digraph in our class then either (1) $D$ has property $Z$ and the reachability relation is not universal; or (2) $D$ does not have property $Z$, the reachability relation is universal and $D$ has infinite in-valency. 
Finally, we describe some distance-transitive, weakly descendant-homogeneous digraphs for which the reachability relation is not universal.
\end{abstract}


\maketitle

\section{Introduction}
A digraph $(D, E)$ consists of a set $D$ of vertices, and a set $E=E(D)$ with $E(D)\subseteq D\times D$ of ordered pairs of vertices, the (directed) edges. We consider digraphs without loops. We will usually refer to the digraph $(D, E)$ simply as `the digraph $D$', and use notation such as `$\alpha \in D$' to indicate that $\alpha$ is a vertex of the digraph $D$.
The \textit{out-valency} of a vertex $\alpha \in D$ is the cardinality of the set $\{u\in D : (\alpha,u)\in E(D)\}$;  similarly, the \textit{in-valency} of $\alpha$ is the cardinality of the set $\{u\in D : (u,\alpha)\in E(D)\}$. A digraph is  \textit{locally finite} if all in- and out-valencies are  finite.

 Let $s\geq0$ be an integer. An \textit{$s$-arc} from  $u$ to  $v$ in a digraph $D$ is a sequence $u_{0}u_{1} \ldots u_{s}$ of $s+1$ vertices such that $u_0 = u$, $u_{s} = v$ and $(u_{i},u_{i+1})\in E(D)$ for $0\leq i<s$ and $u_{i-1}\neq u_{i+1}$ for $0<i<s$. 
We say that a digraph is $s$-\textit{arc transitive} if its automorphism group is transitive on the set of $s$-arcs. A digraph is \textit{highly-arc-transitive}  if it is $s$-arc transitive for every $s \geq 0$. In the finite case, the only connected highly-arc-transitive digraphs are the directed cycles. On the other hand, an infinite highly-arc-transitive digraph does not contain a directed cycle and we have a large list of examples of such digraphs (see \cite{Amato, AT2, AE, Cameronetal, Devos, Evans1,malnic}) constructed using different methods. An easy example of a highly-arc-transitive digraph is the digraph $Z$ which has  the set  $\mathbb{Z}$ of integers as its vertex set  and edges $(i,i+1)$ for every $i \in \mathbb{Z}$. Note that the digraph $Z$ is an infinite  directed tree with in- and out-valency equal to 1. More generally, any infinite regular directed tree is highly-arc-transitive. 

The class of highly-arc-transitive digraphs is investigated by Cameron et al.~\cite{Cameronetal}. The authors construct many examples, pose several questions  and present  approaches to describing such digraphs. One approach is to consider highly-arc-transitive digraphs which have property $Z$. A digraph $D$ is said to have  {\itshape property $Z$} if there is  a homomorphism from $D$ onto $Z$, that is a mapping $f: D \rightarrow \mathbb{Z}$ such that for every edge $(u,v) \in E(D)$ we have $f(v)=f(u)+1$. An infinite directed regular tree is an example of a digraph with this property. There are several other examples of highly-arc-transitive digraphs with property $Z$ and some can be  found in \cite{Amato}, \cite{ AT2}, \cite{ Devos} and \cite{Moller5}. Praeger  \cite{praeger} shows that  a locally finite highly-arc-transitive digraph has property $Z$ whenever the digraph has unequal  in- and out-valency.  In \cite{Cameronetal} the authors ask whether all locally finite highly-arc-transitive digraphs have property $Z$. A negative answer is given by Malni${\check {\rm c}}$ et al.~\cite{malnic} with the construction of infinite highly-arc-transitive digraphs with equal in- and out-valency and without property $Z$. The digraphs constructed by DeVos et al.~\cite{Devos} provide further examples of locally finite highly-arc-transitive digraphs without property $Z$.





\vspace{0.3cm}

  The {\itshape descendant set} of a vertex $\alpha$ in a digraph, denoted by ${\rm desc}(\alpha)$,  is   the set of all vertices $u$ for which there is an $s$-arc from  $\alpha$  to $u$ for some $s\geq 0$. The {\itshape descendant subdigraph} of  $\alpha$ is the induced subdigraph on the set ${\rm desc}(\alpha)$.
Emms and Evans \cite{EvansJo} construct continuum many non-isomorphic countable, primitive, highly-arc-transitive digraphs, all with isomorphic descendant subdigraphs (namely, rooted trees of finite out-valency). This suggests that a classification of such digraphs is out of the question, even under the very strong assumption of high-arc-transitivity. Nevertheless, Peter Neumann (private communication) suggested that a classification of the descendant subdigraphs in these digraphs might be possible. 
 
 In \cite{Amato} we started an investigation of the descendant subdigraphs in highly-arc-transitive digraphs of finite out-valency. We isolate a small number of quite simple properties (P0, P1 and P2 below) satisfied by such descendant subdigraphs, and show that these properties entail rather strong structural consequences. In particular, the descendant subdigraph admits a non-trivial, finite-to-one homomorphism onto a tree.
 
 Let $\Gamma$ be a digraph and let $u$ be a vertex in $\Gamma$. Denote by $\Gamma(u)$ the descendant subdigraph of $u$ and  for $i \geq 0$,  let $\Gamma^i(u)$ denote the set of vertices in $\Gamma$ which can be reached by an $i$-arc from $u$.  Properties P0, P1 and P2 are as follows: 
 \begin{enumerate}
\item[\textbf{P0}]$\Gamma = \Gamma(\alpha)$ is a rooted (root $\alpha$) digraph  with finite out-valency $m > 0$ and $\Gamma^s(\alpha) \cap \Gamma^t(\alpha) = \varnothing$ whenever $s\neq t$.

\item[\textbf{P1}]$\Gamma(u)$ is isomorphic to $\Gamma $ for all vertices $u$ in $\Gamma$.

\item[\textbf{P2}]For $i \in \mathbb{N}$ the automorphism group ${\rm Aut}(\Gamma)$ is transitive on $\Gamma^i(\alpha)$.
\end{enumerate}
 Explicit examples of digraphs having the given properties, but which are not trees, are constructed in  \cite{Amato, Moller4}; and (imprimitive) highly-arc-transitive digraphs having these as descendant subdigraphs are constructed in  \cite{Amato, AT1, AT3}. Moreover, it is shown in  Corollary 4.4 of {\cite{AT3}} that if $\Gamma$ is a  digraph having the given properties, then there is a countable highly-arc-transitive digraph with $\Gamma$ as its descendant subdigraph. 

 
 \vspace{0.2cm}
 Following \cite{Lam}, we say that a digraph $D$ is \textit{distance-transitive} if for every $s \geq 0$, the automorphism group of $D$ is transitive on pairs $(u,v)$ for which there is an $s$-arc from $u$ to $v$, but no $t$-arc for $t <s$. Note that this implies vertex and edge transitivity, but is weaker than being highly-arc-transitive. It is then natural to ask whether there exists a distance-transitive digraph that is not highly-arc-transitive. The answer is positive and an example can be constructed as follows. Start with an infinite locally finite distance-transitive undirected graph $X$. The class of such graphs has been classified by Macpherson \cite{dugald}. Replace each (undirected) edge $\{x,y\}$ with two directed edges, $(x,y)$ and $(y,x)$, to obtain a directed graph $\overrightarrow{X}$.
As noted in \cite{Lam}, distance-transitivity of $\overrightarrow{X}$ follows from distance-transitivity of $X$. However, $\overrightarrow{X}$ is not highly-arc-transitive since it is infinite and contains a finite directed cycle.

 In \cite{AE}, we reprove some of the results of \cite{Amato} in the context of distance-transitive digraphs. In Corollary 2.5 of  {\cite{AE}}, we show that the descendant subdigraph in a distance-transitive digraph $D$ of finite out-valency satisfies P0, P1 and P2; and in Corollary 2.15 of {\cite{AE}} we show that there are only countably many digraphs satisfying properties P0, P1 and P2. Together with Corollary 4.4 of {\cite{AT3}}, the results of \cite{AE} give a reasonable picture of the descendant subdigraphs in distance-transitive digraphs of finite out-valency and infinite in-valency: conditions P0, P1, P2 are necessary and sufficient conditions for a digraph to be a descendant subdigraph in such a digraph, and there are only countably many digraphs satisfying these conditions. 
 
 We also considered descendant subdigraphs under the additional assumption of primitivity. In Corollary 2.5 of {\cite{AE}} we show that if a distance-transitive digraph of finite out-valency  is primitive on vertices, then its descendant subdigraph $\Gamma=\Gamma(\alpha)$ satisfies one further property: 
  \begin{enumerate}
 \item[\textbf{P3}]For $i \in \mathbb{N}$, we have $|\Gamma^i(\alpha)|<|\Gamma^{i+1}(\alpha)|$.
  \end{enumerate}
 In particular, we prove here the following result for  distance-transitive digraphs of out-valency equal to $p^2$, where $p$ is a prime.   
 


\begin{theorem} \label{prime2}Let $D$ be an infinite, primitive, distance-transitive digraph of finite out-valency with no directed cycles. If the out-valency is equal to $p^2$, where $p$ is a prime, then the descendant subdigraph of $D$ is a (rooted) tree. \end{theorem}
This extends  Corollary 3.5 of {\cite{Amato2018}}  for the class distance-transitive digraphs.

\vspace{0.3cm}

Motivated by the results of \cite{AE}, 
we investigate here what other results on highly-arc-transitive digraphs can be extended to distance-transitive digraphs. 
Additionally, we explore what can be said about distance-transitive digraphs when we add the assumption of  descendant-homogeneity.
The results in Section 5 of {\cite{Amato2}} may be considered as the starting point of this investigation. In Corollary 5.7 of {\cite{Amato2}}, we show that a distance-transitive, descendant-homogeneous digraph of finite out-valency with no directed cycles is imprimitive, and this extends the result in Theorem 3.3 of {\cite{AT2}}.




The class of descendant-homogeneous digraphs was introduced in \cite{AT1}. We say that a subset of the vertex set of a digraph is {\itshape finitely generated} if it is a union of finitely many descendant sets, and we may also say it is the {\itshape descendant set of a finite subset}. If it is the union of ${\rm desc}(x)$, for $x \in X$, we write it as ${\rm desc}(X)$. We say that a subdigraph is {\itshape finitely generated} if its vertex set is finitely generated.  A digraph $D$ is {\itshape descendant-homogeneous} if it is vertex transitive and any isomorphism between finitely generated subdigraphs extends to an automorphism of $D$.  There is also a weaker notion: we say that $D$ is {\itshape weakly descendant-homogeneous} if it is vertex transitive and for any two finite subsets $X$ and $Y$ of vertices and an isomorphism $f$ from the subdigraph on ${\rm desc}(X)$ to the subdigraph on ${\rm desc}(Y)$, there is an automorphism of $D$ agreeing with $f$ on $X$. Thus, for descendant-homogeneity, the automorphism must agree with $f$ on the whole of ${\rm desc}(X)$, and for weak descendant-homogeneity, just on $X$. Clearly descendant-homogeneity implies weak descendant-homogeneity.

In Sections \ref{imprimitive} and  \ref{classification}, we turn our attention to distance-transitive digraphs which are, additionally, weakly descendant-homogeneous.
As mentioned earlier, Praeger \cite{praeger} gives a sufficient condition for a locally finite highly-arc-transitive digraph to have property $Z$, namely that the digraph has unequal  in- and out-valency.  
 In Section \ref{imprimitive}, we prove Theorem  \ref{locallyfiniteD} which shows  that weak descendant-homogeneity is a sufficient condition for a locally finite distance-transitive digraph to have property $Z$.

\begin{theorem} \label{locallyfiniteD}Let $D$ be a locally finite, distance-transitive digraph with no directed cycles. If $D$ is weakly descendant-homogeneous, then $D$ has property $Z$. 
\end{theorem}

\noindent There are distance-transitive descendant-homogeneous digraphs with property $Z$ which are not locally finite. An example is the digraph $D(m, \aleph_0)$ constructed in Example 4.4 of  {\cite{Amato2}} (Example \ref{example1} here), where $m\geq 1$ is an integer.  


\vspace{0.3cm}

Another approach used in \cite{Cameronetal} to classify highly-arc-transitive digraphs is the  {\itshape reachability relation}.  An \textit{alternating walk} in a digraph $D$ is a sequence $x_0x_1\ldots  x_n$ of vertices in $D$ such that either $(x_{2i-1},x_{2i})$ and $(x_{2i+1},x_{2i})$ are edges for all $i$, or $(x_{2i},x_{2i-1})$ and $(x_{2i},x_{2i+1})$ are edges for all $i$.  If $e$ and $e'$ are edges in $D$ and there is an alternating walk $x_0x_1\ldots  x_n$ such that $e=(x_0,x_1)$ and either $e'=(x_{n-1},x_n)$ or $e'=(x_n,x_{n-1})$, then $e'$ is said to be {\itshape reachable} from $e$ by an alternating walk and this is denoted by $e \mathcal{A} e'$. Clearly $\mathcal{A}$ is an equivalence relation on the edge set of $D$ which is preserved by the automorphism group of $D$. For an edge $e$ in $D$, we denote its equivalence class by $\mathcal{A}(e)$. A subdigraph spanned by one of the equivalence classes  is called an  {\itshape{alternet}}. When $D$ is 1-arc transitive,  all the alternets are isomorphic to some fixed digraph $\Delta(D)$. If the reachability relation has more than one class then from Proposition 1.1 of {\cite{Cameronetal}} it follows that the digraph $\Delta(D)$ is bipartite. 

Property $Z$ is closely related to the  reachability relation. In a digraph with property $Z$, the reachability relation is not universal  (in fact, there are infinitely many reachability classes). 
In general, the converse is not true. Indeed, Malni${\check {\rm c}}$ et al.~\cite{malnic} construct locally finite, highly-arc-transitive digraphs  without property $Z$ for which the reachability relation is not universal. However, as  stated in the following theorem, we show that  the converse holds for the class of distance-transitive, weakly descendant-homogeneous digraphs. Also, for such a class of digraphs, having property $Z$ is equivalent to the descendant subdigraph not satisfying property P3.

\begin{theorem}  \label{equivalence1}Let $D$ be a distance-transitive, weakly descendant-homogeneous digraph of finite out-valency $m>0$ with no directed cycles. Let $\Gamma$ be the descendant subdigraph of $D$. Then the following are equivalent: 
    \begin{enumerate}[a)]
\item  $D$ has property $Z$.
\item  The reachability relation is not universal.
\item $\Gamma$  does not satisfy  P3.
\end{enumerate} 
\end{theorem}

We may summarise our results as follows. 
\begin{corollary}  \label{summary}Let $D$ be a distance-transitive, weakly descendant-homogeneous digraph of finite out-valency $m>0$ with no directed cycles. Let $\Gamma$ be the descendant subdigraph of $D$. Then $D$ is imprimitive and either:
    \begin{enumerate}[a)]
    \item Case  {\rm(I)}: $D$ has property $Z$,  the reachability relation is not universal and  $\Gamma$  does not satisfy P3;  or
\item  Case  {\rm(II)}: $D$ does not have property $Z$,  the reachability relation is universal, $\Gamma$  satisfies  P3 and  $D$ has infinite in-valency.
\end{enumerate} 
\end{corollary}

\noindent As mentioned earlier, the fact that $D$ is imprimitive is proved in Corollary 5.7 of  {\cite{Amato2}}, and we include a proof here (Corollary \ref{prim}). Examples of digraphs in Case (II) were constructed in \cite{AT2,AT3,Evans1} using a Fra\"{i}ss\'{e}-type construction. In particular,  in Section 4.2 of {\cite{AT2}}, the authors construct uncountably many non-isomorphic  highly-arc-transitive (hence distance-transitive) weakly descendant-homogeneous digraphs without property $Z$. This suggests that the classification of digraphs in Case (II) is out of the question. However,  a classification could still be possible under the stronger assumption of descendant-homogeneity. 

On the other hand, the only  examples known to the author of digraphs which lie in Case (I) are the digraphs described here in Example \ref{example1}. Considering this, in Section  \ref{classification} we work towards a classification of the digraphs in Case (I). Let $D$ be a distance-transitive, weakly descendant-homogeneous digraph with non-universal reachability relation.  Then all  alternets of $D$ are isomorphic to some fixed bipartite digraph $\Delta(D)$. We show that if $\Delta(D)$ is isomorphic to $\Delta(D')$, where $D'$ is a known example in Case (I), then $D$ is isomorphic to $D'$ (see Theorem \ref{description}).

\vspace{0,2cm}


The structure of this paper is as follows. In Section \ref{notation} we establish notation beyond what has already been set in this introduction. In Section \ref{prelim} we recall results of \cite{AE} on the structure of descendant subdigraphs in distance-transitive digraphs  and prove Theorem \ref{prime2} (Section \ref{primes}). In Section \ref{imprimitive} we consider distance-transitive digraphs which are weakly descendant-homogeneous and prove Theorems \ref{locallyfiniteD} (Section \ref{teorema1.2}) and \ref{equivalence1} (Section \ref{teorema1.3}). Finally, in  Section  \ref{classification} we investigate the digraphs in Case (I) of Corollary \ref{summary}. 

\section{Notation and Terminology}\label{notation}
Let $(D, E)$ be a digraph. We will think of a  subset $X$ of the set $D$ of vertices as a digraph in its own right by considering the full induced subdigraph on $X$ (so $E(X) = E(D) \cap X^2$).  Throughout this paper, `subdigraph' will mean `full induced subdigraph'. Thus henceforth, we will not usually distinguish notationally between a digraph and its vertex set. In particular, we will usually refer to the digraph $(D, E)$ simply as `the digraph $D$'. 
Furthermore, we will use notation such as `$\alpha \in D$' to indicate that $\alpha$ is a vertex of the digraph $D$.

We denote the automorphism group of a digraph $D$ by ${\rm Aut}(D)$. We say that $D$ is \textit{transitive} (or \textit{vertex transitive})   if   ${\rm Aut}(D)$ is transitive on $D$. We say that $D$ is \textit{edge transitive}   if   ${\rm Aut}(D)$ is transitive on $E(D)$. Note that the set of 1-arcs of $D$ is  the set of edges of $D$, so the notions of edge transitivity and 1-arc transitivity   are the same.

The \textit{out-valency} of a vertex $\alpha \in D$ is the cardinality of the set 
$$
{\rm out}_D(\alpha)=\{u\in D : (\alpha,u)\in E(D)\}
$$ 
of out-vertices (or out-neighbours) of $\alpha$;  similarly, the \textit{in-valency} of $\alpha$ is the cardinality of the set 
$$
{\rm in}_D(v)=\{u\in D : (u,\alpha)\in E(D)\}
$$ 
of in-vertices (or in-neighbours). We write simply ${\rm out}(\alpha)$ and ${\rm in}(\alpha)$  for the above sets when it is clear  the digraph $D$ we are referring to. A digraph with both in- and out-valencies finite is said to be \textit{locally finite}. 
If $D$ is transitive then all out-valencies are equal and we call this the `out-valency of $D$'. Similarly for in-valencies.

A {\itshape path} in a digraph $D$ is a sequence $v_0 v_1 \ldots v_n$ of distinct vertices such that either $(v_i, v_{i-1}) \in E(D)$ or $(v_{i-1}, v_i) \in E(D)$ for all $i$. A digraph is said to be {\itshape connected}  if for any pair $u,v$ of vertices there is a path from $u$ to $v$. 

For an integer $s\geq0$, an \textit{$s$-arc} from  $u$ to  $v$ in a digraph $D$ is a sequence $u_{0}u_{1} \ldots u_{s}$ of $s+1$ vertices such that $u_0 = u$, $u_{s} = v$ and $(u_{i},u_{i+1})\in E(D)$ for $0\leq i<s$ and $u_{i-1}\neq u_{i+1}$ for $0<i<s$. An $s$-arc is also called a directed path of length $s$.  A sequence $u_{0}u_{1} u_2 \ldots $ such  that $(u_{i},u_{i+1})\in E(D)$ for $i \geq 0$ is called an \textit{infinite arc}. A sequence $\ldots u_{-2}u_{-1}u_{0}u_{1} u_2 \ldots $ such  that $(u_{i},u_{i+1})\in E(D)$ for $i \in \mathbb{Z}$ is called a \textit{two-way infinite arc}. Two-way infinite arcs are often called  \textit{directed lines}. 

Let $s\geq0$ be an integer and let $u \in D$. We denote by $D^s(u)$ (or ${\rm desc}^s(u)$) the set of vertices of $D$ which are reachable by an $s$-arc from $u$. The \textit{descendant set} $D(u)$ (or ${\rm desc}(u)$) of $u$ is $\bigcup_{s\geq 0} D^s(u)$.  Similarly, the set ${\rm anc}(u)$ of {\textit {ancestors}} of $u$ is the set of vertices of which $u$ is a descendant, and $D^{-s}(u)$ (or    ${\rm anc}^s(u)$) is the set of vertices $x$ in $D$ for which there is an $s$-arc from $x$ to $u$. For a subset $X\subseteq D$, we write $D^s(X)=\bigcup_{x\in X} D^s(x)$ and $D(X)=\bigcup_{x\in X} D(x)$, and refer to the latter as the descendant set of $X$.  The {\itshape descendant subdigraph} of $u$  is the induced subdigraph on  the set $D(u)$.

 Let $\alpha \in D$, and let $\Gamma$ be the descendant subdigraph  of $\alpha$. If ${\rm Aut}(D)$ is transitive on the set of vertices of $D$, then the  subdigraph $D(u)$ is isomorphic to $\Gamma $ for all vertices $u\in D$, and we shall speak of the digraph $\Gamma$ as \textit{the descendant subdigraph} of $D$.
 We will be considering this as a digraph with its full induced structure from $D$. 
 We refer to $\alpha$ as the \textit{root} of $\Gamma$ and write $\Gamma = \Gamma(\alpha)$ to indicate that any vertex of $\Gamma$ is a descendant of $\alpha$. Similarly, we write $\Gamma^s$, or  $\Gamma^s(\alpha)$, for the set of vertices reachable by an $s$-arc starting at $\alpha$ and  if $\beta \in \Gamma(\alpha)$, then  we write $\Gamma(\beta)$ for the subdigraph of $\Gamma$ induced on $ {\rm desc}(\beta) \subseteq \Gamma(\alpha)$. It is clear that if $D$ is highly-arc-transitive, then ${\rm Aut}(\Gamma(\alpha))$ is transitive on $s$-arcs in $\Gamma(\alpha)$ which start at $\alpha$. Similarly, if $D$ is distance-transitive, then ${\rm Aut}(\Gamma(\alpha))$ is transitive on $\Gamma^n(\alpha)$ for each $n \in \mathbb{N}$.

Let $G$ be a group acting  transitively on a set $\Omega$. The image of an element $a \in \Omega$ under $g \in G$ will be denoted by $a^g$. The {\itshape stabilizer} of $a \in \Omega$ is the subgroup $G_a=\{g \in G \mid a^g=a\}$. For a subset $X \subseteq \Omega$, the {\itshape setwise stabilizer} of $X$ is the subgroup 
$G_{\{X\}}=\{g \in G \mid X^g=X\}$, and the {\itshape pointwise stabilizer} is the subgroup $G_{(X)}=\{g \in G \mid x^g=x {\mbox { for all }} x \in X\}$.  For a subset  $Y$  of $\Omega$ which is invariant under the action of $G$, we denote by $G^Y$ the {\itshape permutation group induced by the action of $G$ on the set $Y$}.
A proper subset $B$ of $\Omega$ with at least two elements is called a {\itshape block of imprimitivity} for $G$  if for every $g \in G$ either $B^g=B$ or $B \cap B^g=\varnothing$. An equivalence relation on $\Omega$ that is preserved by $G$ is called a $G$-{\itshape congruence}. If $B$ is a block of imprimitivity then $B$ and its translates, $B^g$ with $g \in G$, are the classes of a $G$-congruence. Conversely, if we have a non-trivial proper $G$-congruence then each one of its classes is a block of imprimitivity. If there are no blocks of imprimitivity or, equivalently, there are no non-trivial proper $G$-congruences, we say that the action of $G$ on $\Omega$ is {\itshape primitive}. 

We say that a digraph $D$ is \textit{primitive} if ${\rm Aut}(D)$ is primitive on the vertex set of $D$.  Otherwise, we say $D$ is \textit{imprimitive}.


For a digraph $D$ and an equivalence relation $\rho$ on the vertex set of $D$, the \textit{quotient digraph} $D/ \rho$ is the digraph which has the set of $\rho$-classes as the vertex set, and if $A$ and $B$ are two $\rho$-classes then $(A,B)$ is an edge in $D/ \rho$ if, and only if, there exist vertices $a \in A$ and  $b \in B$ such that $(a,b)$ is an edge in $D$.

Throughout the paper our digraphs are always assumed to be infinite and connected. 

\section{Distance-transitive digraphs of finite out-valency}\label{prelim}

In \cite{Amato} we isolate a small number of simple properties, P0, P1, and P2, satisfied by the descendant subdigraphs in highly-arc-transitive digraphs of finite out-valency, and show that these properties entail rather strong structural consequences. Later, in \cite{AE},  we reprove some of the results of \cite{Amato} in the context of  distance-transitive digraphs. 
In this section, we recall some of the results in Section 2.1 of {\cite{AE}} on the structure of descendant subdigraphs in distance-transitive digraphs of finite out-valency, and prove Theorem \ref{prime2}.


\subsection{Preliminaries}\label{prelim1}


 For a digraph  $\Gamma$, consider the following properties:

\begin{enumerate}
\item[\textbf{P0}]$\Gamma = \Gamma(\alpha)$ is a rooted digraph with finite out-valency $m > 0$ and $\Gamma^s(\alpha) \cap \Gamma^t(\alpha) = \varnothing$ whenever $s\neq t$.

\item[\textbf{P1}]$\Gamma(u)$ is isomorphic to $\Gamma $ for all $u\in \Gamma$.

\item[\textbf{P2}]For $i \in \mathbb{N}$ the automorphism group ${\rm Aut}(\Gamma)$ is transitive on $\Gamma^i(\alpha)$.

\item[\textbf{P3}]For $i \in \mathbb{N}$, we have $|\Gamma^i(\alpha)|<|\Gamma^{i+1}(\alpha)|$.
\end{enumerate}

Properties P0, P1, P3 are the same as G0, G1, G2 in Section 2.1 of {\cite{AE}}, respectively.  In the latter, the authors also consider the following property: 
\begin{enumerate}
\item [\textbf{G3}] There is an integer $k \geq 1$ such that if $l \geq k$ and $x \in \Gamma^l(\alpha)$ and $z \in \Gamma(x)$, then ${\rm anc}(z) \cap \Gamma^1(\alpha)={\rm anc}(x) \cap \Gamma^1(\alpha)$.
\end{enumerate}

\begin{lemma}\label{G3}
 ({\cite[Lemma 2.3]{AE}}) Suppose $\Gamma$ satisfies P0, P1 and P2. Then $\Gamma$ satisfies G3. 
\end{lemma}

A {\itshape{priori} }there could be continuum-many isomorphism types of digraphs with the properties above. Theorem 2.15 of {\cite{AE}} shows that there are only countably many isomorphism types of digraphs $\Gamma$ satisfying P0, P1 and G3. To establish this, the authors show that there is a natural equivalence relation $\rho$ on $\Gamma$ (refining the  `layering' of $\Gamma$ given by P0) such that the quotient digraph $\Gamma / \rho$ is a directed tree, as we shall see in Section \ref{structure}. If P3 also holds then $\Gamma / \rho$ is not a directed line and the size of the layers $\Gamma^n(\alpha)$ grows exponentially.

The class of digraphs $\Gamma$ satisfying P0, P1 and G3 includes the descendant subdigraphs of some distance-transitive digraphs, as the next result shows. Note that property P2 follows directly from distance-transitivity. 
\vspace{0.3cm}


\begin{corollary}\label{desc}
 ({\cite[Corollary 2.5]{AE}}) Suppose $D$ is a distance-transitive digraph of finite out-valency $m>0$ and is either of infinite in-valency, or has no directed cycles. Then the descendant subdigraph $\Gamma$ in $D$ satisfies P0, P1, P2 and G3. If the automorphism group of $D$ is also primitive on vertices, then $m>1$ and $\Gamma$ satisfies P3.
\end{corollary}


The result shows that P0, P1, P2 and P3 are necessary conditions for a digraph to be the descendant subdigraph in a distance-transitive digraph of finite out-valency acting primitively on vertices. However, the properties are not sufficient. Indeed, Evans \cite{Evans1} constructs an infinite, imprimitive, highly-arc-transitive digraph of finite out-valency whose descendant set is a (rooted) regular tree of out-valency 2. Note that a rooted regular directed tree of out-valency greater than 1  satisfies P0 to P3.

In  Theorem 1.2 of {\cite{AE}}, the authors construct a primitive distance-transitive digraph $D$ of finite out-valency and infinite in-valency. By Corollary \ref{desc} the descendant subdigraph of $D$ satisfies P3. 
On the other hand, the digraphs constructed in Example 4.4 of {\cite{Amato2} (Example \ref{example1} here) are  imprimitive distance-transitive digraphs with descendant subdigraph not satisfying P3. 

\vspace{0.3cm}

We end this section with the following definition. 

\begin{definition}\label{delta}
Let $D$ be a digraph and let $n$ be a natural number with $n \geq 1$. For vertices  $u,v$ in $D$, write $\delta_n(u,v)$ if, and only if, 
$$
D^n(u)=D^n(v).
$$
Otherwise, say that $\delta_n(u,v)$ does not hold.
\end{definition}

Clearly $\delta_n$ is an equivalence relation on the vertex set of $D$ which is preserved by ${\rm Aut}(D)$. 
Note that if there are distinct vertices $u,v$ in $D$ such that $D^n(u)=D^n(v)$, then  $D^{l}(u)=D^{l}(v)$ for $l \geq n$. This shows that  if $\delta_n$ is non-trivial for some $n$, then $\delta_l$ is non-trivial for all $l \geq n$; and if $\delta_n$ is trivial for some $n$, then $\delta_l$ is trivial for all $l \leq n$. Moreover, if there exists $n$ such that  $\delta_n$ is universal then $D$ contains a directed cycle. Indeed, for $a \in D$ and $v \in D^n(a)$, we have $v \in D^n(v)$ since $D^n(v)=D^n(a)$. This means that there is a directed cycle of length $n$ containing $v$.






\subsection{The structure of descendant subdigraphs}\label{structure}
Throughout this section we assume that $\Gamma$ is a digraph satisfying P0, P1 and  G3. We let $k$ be an integer satisfying the condition in G3. We refer to Section 2.2 of {\cite{AE}} for proofs.

\begin{definition}{\cite[Definition 2.7]{AE}}\label{rho}
Suppose $s$ is a non-negative integer. 
\begin{enumerate}[a)]
\item  Suppose $\beta \in \Gamma, x \in \Gamma^n(\beta)$ and $s\leq n$. Define
$$
\Gamma^{-s}_{\beta}=\{ w \in \Gamma^{n-s}(\beta): x \in {\Gamma}(w)\}.
$$
\item For $l\geq k$ and $x,y \in \Gamma^{l}(\alpha)$ write $\rho(x,y)$ if, and only if, 
$$
\Gamma^{-k+1}_{\alpha}(x)=\Gamma^{-k+1}_{\alpha}(y).
$$
(Say that $\rho(x,y)$ does not hold in all other cases.)
\end{enumerate}
\end{definition}
So for $x,y \in \Gamma^l(\alpha)$ we have that $\rho(x,y)$ holds if, and only if, $x,y$ have the same ancestors in $\Gamma^{l-k+1}(\alpha)$. Clearly $\rho$ is an  equivalence relation on $\cup_{l\geq k}\Gamma^l(\alpha)$ which is invariant under the action of ${\rm Aut}(\Gamma)$.

For $l\geq k$ and $x \in \Gamma^l(\alpha)$ we write $[x]_{\rho}$ for the $\rho$-equivalence class containing $x$. We use notation such as  ${\bf{v}},{\bf{w}}$ etc. for such classes and write $\Gamma( {\bf{u}})=\cup_{x \in {\bf{u}}} {\Gamma}(x)$ and $\Gamma^s( {\bf{u}})=\cup_{x \in {\bf{u}}} {\Gamma}^s(x)$. 

\begin{lemma}{\cite[Lemma 2.9]{AE}}\label{2.9}
Suppose $l \geq k$ and ${\bf{v}} \subseteq \Gamma^l(\alpha)$ is a $\rho$-class. Let $w \in \Gamma({\bf{v}})$. Then $[w]_{\rho}\subseteq \Gamma( {\bf{v}})$.
\end{lemma}

\begin{corollary} {\cite[Corollary 2.10]{AE}}\label{2.10}
Suppose $l \geq k$ and $v \in  \Gamma^l(\alpha)$. Let ${\bf{v}}$ be the $\rho$-class containing $v$. Then the quotient digraph $\Gamma( {\bf{v}}) / \rho$ is a rooted directed tree with finite out-valencies. 
\end{corollary}

Note that for $\beta \in \Gamma(\alpha)$ we can consider the equivalence relation $\rho$ computed in both $\Gamma(\alpha)$ and $\Gamma(\beta)$, where in the latter we only consider ancestors in  $\Gamma(\beta)$ when defining $\rho$: a {\itshape{priori}} this gives a coarser relation.

\begin{lemma}{\cite[Lemma 2.11]{AE}}
Suppose $\beta \in \Gamma^n(\alpha)$ and $x \in \Gamma^l(\beta)$ with $l \geq 2k-1$. Then the $\rho$-class containing $x$ is the same whether it is computed in $\Gamma(\alpha)$ or $\Gamma(\beta)$.
\end{lemma}

Let $l\geq 2k-1$ and let ${\bf{v}}$ be a $\rho$-class in $\Gamma^l(\alpha)$. Let $T({\bf{v}})$ be the structure consisting of the induced digraph on $\Gamma({\bf{v}})$ together with the equivalence relation induced by $\rho$ (coming from $\Gamma(\alpha)$). Recall that by Lemma \ref{2.9},  $T({\bf{v}})$ is a union of $\rho$-classes in $\Gamma(\alpha)$. If  ${\bf{w}}$ is another $\rho$-class (in $\cup_{l\geq 2k-1}\Gamma^l(\alpha)$) then by a $\rho$-{\itshape{isomorphism}} between $T({\bf{v}})$ and $T({\bf{w}})$ we mean a digraph isomorphism which respects $\rho$.

\begin{corollary} {\cite[Corollary 2.12]{AE}}\label{2.12}
Suppose ${\bf{v}}$ is a $\rho$-class in $ \Gamma^l(\alpha)$ with $l \geq 2k-1$. Then there is a $\rho$-class ${\bf{w}}$ in $\Gamma^{2k-1}(\alpha)$ and a $\rho$-isomorphism from $T({\bf{w}})$ to $T({\bf{v}})$.
\end{corollary}

Thus to any digraph $\Gamma$ satisfying P0, P1, G3, there are associated a finite number of $\rho$-isomorphism types of $T({\bf{v}})$. In particular, we can refine Corollary \ref{2.10} to: 

\begin{corollary} {\cite[Corollary 2.13]{AE}}\label{2.13}
Suppose $l \geq 2k-1$ and  ${\bf{v}}\subseteq \Gamma^l(\alpha)$ is a $\rho$-class.  Then the quotient digraph $T( {\bf{v}}) / \rho$ is a rooted directed tree with a  finite number of out-valencies. 
\end{corollary}

If we now assume that   $\Gamma$ also satisfies P2, then we  have the following result. 
\begin{theorem}\label{regular}
Let $\Gamma$ be a digraph satisfying P0, P1 and P2. Then there is an integer $k$ such that for $l \geq 2k-1$ and ${\bf{w}} \in \cup_{l\geq 2k-1}\Gamma^l(\alpha)$, the quotient digraph $T( {\bf{w}}) / \rho$ is a rooted directed  tree of out-valency  $s \geq 1$. In particular, the subdigraph induced on ${\rm desc}(\Gamma^l(\alpha))$ is the disjoint union of finitely many copies of a fixed digraph $T$, where $T$ is isomorphic to $T( {\bf{w}})$.
\end{theorem}
\begin{proof}
By Lemma \ref{G3}, $\Gamma$ satisfies G3. Let $k$ be the integer satisfying the condition in G3, and let  ${\bf{w}}$ be a $\rho$-class in $\Gamma^{2k-1}(\alpha)$. By property  P2 and Corollary \ref{2.12}, for $l \geq 2k-1$ and ${\bf{v}}$ a $\rho$-class in $\Gamma^{l}(\alpha)$ there is a $\rho$-isomorphism from $T( {\bf{w}})$ to $T( {\bf{v}})$. The result now follows from Corollary \ref{2.13}.
\end{proof}




\subsection{Property P3}\label{notP3}

Let $\Gamma$ be a digraph satisfying P0, P1 and P2. We prove here that, for $i \geq 1$, the size of the layer  $\Gamma^i$ is smaller than or equal to the size of the layer $\Gamma^{i+1}$. So there are only two cases to consider: that   $\Gamma$ satisfies P3  and that  $\Gamma$ does not satisfy P3. 
 We then  investigate the structure of $\Gamma$ in each case.

Let $i \geq 1$ and $u\in \Gamma$. By property P2, the in-valency of any two vertices in $\Gamma^i$ is the same, and we denote this in-valency by $r_i$. By property P1,  the in-valency of vertices in $\Gamma^i(u)$ within ${\rm desc}(u)$ is equal to $r_i$. Let $u \in \Gamma^1$. Then $\Gamma^{i}(u) \subseteq \Gamma^{i+1}$ and therefore $r_i \leq r_{i+1}$. So
\begin{equation}\label{sequence}
1=r_1\leq r_2 \leq \ldots \leq r_n\leq \ldots
\end{equation}
 is an infinite non-decreasing sequence of natural numbers, which we refer to as the `in-valency sequence of $\Gamma$'. Also by P1, we have $|\Gamma^i|=|\Gamma^i(u)|$. Since the set $\Gamma^{i}(u)$ is contained in $ \Gamma^{i+1}$, we conclude the following: 
 
 \begin{lemma} ({\cite[Lemma 3.4]{Amato}})\label{3.4}
For $i \in \mathbb{N}$, we have $|\Gamma^i|\leq |\Gamma^{i+1}|$. 
  \end{lemma}

 \begin{lemma}\label{ri}
  For $i \geq1$, $r_{i+1} \leq m$ and equality holds if, and only if, $|\Gamma^{i}|=|\Gamma^{i+1}|$.
 \end{lemma}
 \begin{proof}
 By counting edges from $\Gamma^{i}$ to  $\Gamma^{i+1}$ in two different ways, we obtain $|\Gamma^{i}|\times m=|\Gamma^{i+1}| \times r_{i+1}$. The result now follows from Lemma \ref{3.4}.
 \end{proof}
 \vspace{0,2cm}
 
Then (\ref{sequence}) is an infinite non-decreasing sequence of natural numbers smaller than or equal to $m$.
Hence there is a natural number $n \geq 1$ such that $r_j=r_n$ for all $j \geq n$. This allows us to make the following definition. 

\begin{definition}\label{defN}
The natural number  $N=N(\Gamma)$ is the smallest number such that $r_j=r_N$ for all $j\geq N$.  
\end{definition}
Throughout the paper, $N$ will be used for the number as given in this definition.  We refer to $r_N$ as the `{\itshape ultimate in-valency}' of $\Gamma$.

\begin{lemma} ({\cite[Corollary 4.10]{Amato}}) \label{ultimate}
The ultimate in-valency $r_N$ divides the out-valency $m$ of $\Gamma$.
\end{lemma}


\begin{remark} 
If  $\Gamma$ is a digraph satisfying P0 to  P3 then its out-valency $m$ is greater than 1. Indeed, if $m=1$ then $|\Gamma^i|=1$ for all $i$,  contradicting P3. 
\end{remark}

\begin{lemma}\label{invalency}
({\cite[Lemma 3.3]{Amato2}}) Let $\Gamma$ be a digraph satisfying P0, P1 and P2. Suppose $\Gamma$ does not satisfy P3. Then 
\begin{enumerate}[a)]
\item  for  $i \geq N-1$, $|\Gamma^i|=|\Gamma^{N-1}|$;
\item   for $i \geq N$, $r_i=m$; and 
\item  for $t \geq 1$, $u,v \in \Gamma^{t}$ and $l \geq N-1$, we have $\Gamma^l(u)=\Gamma^{l+t}=\Gamma^l(v)$.
\end{enumerate}
\end{lemma}
\begin{proof}
Since $\Gamma$ does not satisfy  P3, there is a smallest number $l \geq 1$ such that $|\Gamma^{l-1}|=|\Gamma^{l}|$. Then, by Lemma \ref{ri},   $r_l=m$. This, together with the fact that $r_i\leq r_{i+1}$ and $r_i\leq m$  for all $i$, implies that  $r_i=m$ for all $i \geq l$. Now, by the choice of $l$ and Lemma \ref{ri}, we have  $r_{l-1}<m$. Hence $l$ is the least such that $r_j=r_l$ for all $j \geq l$ and therefore $l=N$. This proves (a) and (b).

For (c), let  $t \geq 1$, $u,v \in \Gamma^{t}$ and $l \geq N-1$. The sets $\Gamma^l(u)$ and $\Gamma^l(v)$ are both contained in $\Gamma^{l+t}$. By (a),  $|\Gamma^i|=|\Gamma^{N-1}|$ for $i \geq N-1$. So $|\Gamma^{l+t}|=|\Gamma^l|$. Now, by P1, the sets $\Gamma^l(v),\Gamma^l(u)$ and $\Gamma^l$ have the same cardinality. Hence $\Gamma^l(u)=\Gamma^{l+t}=\Gamma^l(v)$, since  both  $\Gamma^l(u)$ and $\Gamma^l(v)$ are contained in $\Gamma^{l+t}$.
\end{proof}



\begin{proposition} ({\cite[Proposition 3.5]{Amato2}})\label{caracterizaP3}
Let $\Gamma$ be a digraph satisfying P0, P1 and P2. Then $\Gamma$ satisfies P3  if, and only if, there are vertices $x,y \in \Gamma$ such that ${\rm desc}(x) \cap {\rm desc}(y) =\varnothing$.
\end{proposition}



The next result appears in Lemma 4.1 of  {\cite{Amato2}} in the context of highly-arc-transitive digraphs. The proof depends only on Lemma \ref{invalency} and on Proposition \ref{caracterizaP3}, and both hold for digraphs $\Gamma$ satisfying P0, P1 and P2. Since by Corollary \ref{desc} the descendant subdigraph in a distance-transitive digraph of finite out-valency satisfies properties P0 to P2, we may state the following.



\begin{lemma}\label{41} ({\cite[Lemma 4.1]{Amato2}}) Let $D$ be a distance-transitive digraph of finite out-valency $m>0$ with no directed cycles such that its descendant subdigraph does not satisfy  P3. Let $L$ be a directed line in $D$ and denote by $F$ the subdigraph induced on the descendant set of $L$. Then $F$ has property $Z$, ${\rm out}_F(v)={\rm out}_D(v)$  and $ |{\rm in}_F(v)|=m$ for all $v \in F$. Moreover, $F$ has exactly two ends.
\end{lemma}


Let $X$ be an infinite undirected graph, and let $P$ and $Q$ be two one-way infinite paths (also called rays) in $X$. $P$ and $Q$ are \textit{equivalent} in $X$ if there are infinitely many finite disjoint paths connecting vertices in $P$ to vertices in $Q$. The equivalence classes of all infinite paths with respect to this relation are called the \textit{ends} of $X$.  It is well known that every infinite vertex transitive graph has 1, 2, or infinitely many ends. See  \cite{diestel2} for more on ends of graphs.

For an infinite digraph $D$, an \textit{end} of $D$ is an end of the underlying undirected graph of $D$ (that is, the graph obtained from $D$ by replacing each directed edge $(a,b)$ with an undirected edge $\{a,b\}$).  The digraph $Z$ is an example of a digraph with two ends, and an infinite directed tree of out-valency (or in-valency)   greater than 1 is an example of a digraph with infinitely many ends. 


\begin{corollary} \label{cor41}Let $D$ be a distance-transitive digraph of finite out-valency $m>0$ with no directed cycles. Suppose that its descendant subdigraph does not satisfy  P3. Then the in-valency of $D$ is greater than or equal to $m$. In particular, if the in-valency is equal to $m$ then $D$ is the digraph induced on the descendant set of a directed line.
\end{corollary}
\begin{proof}
Let $t$ be the in-valency of $D$. Let $L$ be a directed line in $D$, and denote by $F$ the subdigraph induced on  the descendant set of $L$.  Since ${\rm in}_F(v)$ is subset of ${\rm in}_D(v)$ and, from Lemma \ref{41}, we have $|{\rm in}_F(v)|=m$, it follows that $t\geq m$. If equality holds, then ${\rm in}_F(v)={\rm in}_D(v)$ and we conclude that  $D=F$. 
\end{proof}

\vspace{0.3cm}

Praeger \cite{praeger}  shows that if $D$ is an infinite, connected, vertex  and edge transitive digraph with finite but unequal in- and out-valency then $D$ has property $Z$. We use this fact and Corollary \ref{cor41} to prove the following. 

\begin{corollary}\label{finitevalency}Let $D$ be a distance-transitive digraph  of finite out-valency $m>0$ with no directed cycles. Suppose that  its descendant subdigraph does not satisfy  P3. If $D$  has finite in-valency then $D$ has property $Z$.
\end{corollary}
\begin{proof} Let $t$ denote the in-valency of $D$. From Corollary \ref{cor41} we know that $t \geq m$, and that if $t=m$ then $D$ has property $Z$ (since $D$ is the digraph induced on the descendant set of a directed line). Now assume $t>m$. Then $D$ is a locally finite edge transitive digraph with finite but unequal in- and out-valency. It then follows from \cite{praeger} that $D$ has property $Z$. 
\end{proof}

\vspace{0,3cm}

We end this section with a result on  two-ended highly-arc-transitive digraphs. The general structure of two-ended transitive graphs has been studied in several papers \cite{ Devos, diestel, Moller5}. 
The following result may be a direct consequence of the results in \cite{ Devos, diestel, Moller5} but as we have not found an explicit proof, we include one here.  

  \begin{theorem} \label{2end} Let $D$ be a two-ended highly-arc-transitive digraph. Then  $D$ is the digraph induced on  the descendant set of a directed line in $D$.  
   \end{theorem}
 \begin{proof} Let $D$ be as stated and let $X$ be the underlying undirected graph of $D$. Let $A={\rm Aut}(X)$. Then $X$  is a two-ended graph and $A$ acts transitively on the set of vertices of $X$. Then, by  \cite[Theorem 7]{diestel},  $X$ is locally finite  and therefore $D$ is locally finite. Now, by   \cite[Lemma 7]{Moller5},  there is a homomorphism $\phi$ from $D$ onto $Z$ with $\phi^{-1}(i)$ finite. Moreover, by \cite[Corollary 9]{Moller5}, the in-valency of $D$ is equal to its out-valency, which we will denote by $m$. 
 
 Let $L= \ldots v_{-1} v_0v_1v_2  \ldots$ be a directed line in $D$ such that $v_0 \in  \phi^{-1}(0)$. Now let $\Gamma={\rm desc}(v_0)$. For $j \geq 0$, we have $\Gamma^j \subseteq \phi^{-1}(j)$ so that 
 $|\Gamma^j| \leq |\phi^{-1}(j)|$. Since $\phi^{-1}(j)$ is finite and  $|\phi^{-1}(j)|= |\phi^{-1}(i)|$ for all $i,j$, the  subdigraph $\Gamma$ does not satisfy P3. It then follows from Corollary \ref{cor41} that the subdigraph $F$ induced on the descendant set of the line $L$ has in-valency equal to $m$. Since $D$ has in-valency $m$, we conclude that $F=D$.
\end{proof}


 \subsection{Digraphs satisfying P0 to P3}
 
 Let $\Gamma$ be a digraph satisfying P0 to P3, and let $A={\rm Aut}(\Gamma)$. The goal of this section is to prove that there is a system $\Omega=\{\omega_1,\ldots, \omega_s\}$ of blocks of imprimitivity for $A$ in  $\Gamma^1$, the first layer of $\Gamma$,  with $s\geq 2$ such that for $i \neq j$ the descendant sets  ${\rm desc}(\omega_i)$ and ${\rm desc}(\omega_j)$ do not intersect. In particular, this shows that the digraph $\Gamma \setminus \{\alpha\}$,  obtained from $\Gamma$ by deleting the root $\alpha$, has at least $s$ infinite components. 
 
\begin{proposition}(\cite[Proposition 5.4]{Amato2})\label{equivalence}
Let  $\Gamma$ be a digraph satisfying  P0, P1 and P2. Then the following are equivalent: 
{\rm (a)} $\Gamma$ satisfies P3; {\rm (b)}   there are vertices $x,y$ in $\Gamma$ such that ${\rm desc}(x)\cap {\rm desc}(y)=\varnothing$; {\rm(c)} the quotient digraph  $\Gamma/\rho$ is not a directed line;
{\rm (d)}  $\Gamma$ has more than one end.
\end{proposition}
 

\begin{definition}\label{condition}
For a vertex $x$ in $\Gamma$, $\mathcal{C}_{\Gamma}(x)$ is the following condition: 

\center There are non-empty subsets $U$ and $V$ of $\Gamma^1(x)$ such that ${\rm desc}(U) \cap {\rm desc}(V)=\varnothing$.
\end{definition}

 
 
 
 
 \begin{proposition} \label{conditionstar}Let $\Gamma$ be a digraph satisfying  P0, P1 and P2. Then $\Gamma$  satisfies  P3 if, and only if, $\mathcal{C}_{\Gamma}(x)$ holds for some $x$ in $\Gamma$.
\end{proposition}
 \begin{proof} 
 By Lemma \ref{G3}, $\Gamma$ satisfies G3. Let $k$ be the integer satisfying the condition in G3 and let $l  \geq 2k-1$. Assume that $\Gamma$ satisfies P3. Let $\rho$ be as in Definition \ref{rho}. Now let $\bf{w}$ be a $ \rho$-class in $ \Gamma^l$, and let $T({\bf{w}})$ be the structure consisting of the induced digraph on $\Gamma({\bf{w}})$ together with the equivalence relation induced by $\rho$. By Theorem  \ref{regular}, the quotient digraph  $T(\bf{w})/  \rho$ is a rooted directed regular tree of out-valency  $n\geq 1$. Also, by  Proposition \ref{equivalence},  the quotient digraph $\Gamma/ \rho$ is not a direct line and therefore  $n  \geq 2$. So for every $ \rho$-class $\bf{u}$ in $\bigcup_{i\geq l}\Gamma^i$ there are at least two distinct $ \rho$-classes in $ \Gamma^1(\bf{u})$.
 
 Suppose there is a vertex $x$ in $\bigcup_{i\geq 2k-1}\Gamma^i$ such that the set $\Gamma^1(x)$ intersects (at least) two distinct $\rho$-classes, say $\bf{u }$ and $\bf{v}$.  Then $\mathcal{C}_{\Gamma}(x)$ holds, with $U={\bf{u }} \cap \Gamma^1(x)$ and $V={\bf{v}}\cap \Gamma^1(x)$, since the descendant of any two distinct $\rho$-classes do not intersect.
 
 Now suppose that for every vertex $x$ in $\bigcup_{i\geq 2k-1}\Gamma^i$  the set $\Gamma^1(x)$ is completely contained in some $\rho$-class. Define a relation $R$ on $\Gamma^{2k-1}$ as follows: for vertices $a,b$ write $R(a,b)$ if, and only if,  $\Gamma^1(a)$ and $\Gamma^1(b)$ are contained in a same $\rho$-class  of the level $\Gamma^{2k}$. Then $R$ is an ${\rm Aut}(\Gamma)$-invariant equivalence relation which is coarser than $\rho$ since distinct $\rho$-classes have disjoint descendant sets. Let  $\bf{y}$ be a  $\rho$-class in  $\Gamma^{2k-1}$. We have $n$  distinct $\rho$-classes in $\Gamma^1(\bf{y})$ and each one of them  corresponds to a distinct $R$-class in $\bf{y}$. So there are $n$ distinct $R$-classes in $\bf{y}$. Moreover, since the descendant sets of any two distinct $\rho$-classes in $\Gamma(\bf{y})$ have disjoint intersection, the descendant sets of any two distinct $R$-classes in $\bf{y}$ will also have disjoint intersection. Hence level $\Gamma^{2k-1}$ is a disjoint union of $R$-classes with the property that the descendant sets of  any two distinct such classes have disjoint intersection.  
 
 We now consider the previous level, $\Gamma^{2k-2}$,  and assume that there is $x\in  \Gamma^{2k-2}$ such that $\Gamma^1(x)$ intersects at least two distinct $R$-classes in $\Gamma^{2k-1}$. Then, as above,  $\mathcal{C}_{\Gamma}(x)$ holds. Otherwise, $R$ is an ${\rm Aut}(\Gamma)$-congruence on  $\Gamma^{2k-2}$ which has at least two distinct classes with the property that  any two such classes have disjoint descendant sets. Continuing in this way, we find a vertex $x \in \Gamma^i$, for some $i \in \{0,1, \ldots, 2k-1\}$, for which  $\mathcal{C}_{\Gamma}(x)$ holds.
 
 For the converse, assume there is a vertex $x \in \Gamma$ such that  $\mathcal{C}_{\Gamma}(x)$ holds. Then there are non-empty subsets $U,V$ of $\Gamma^1(x)$ with ${\rm desc}(U) \cap {\rm desc}(V)=\varnothing$. For $a \in U$ and $b \in V$, we have ${\rm desc}(a) \cap {\rm desc}(b)=\varnothing$ and the result  follows from Proposition \ref{equivalence}.
 \end{proof}

  \begin{corollary} \label{conditionstar2}Let $D$ be a distance-transitive digraph of finite out-valency  $m>0$ with no directed cycles. The descendant subdigraph $\Gamma$  of $D$ satisfies  P3 if, and only if, $\mathcal{C}_{\Gamma}(x)$ holds for some $x \in \Gamma$.
\end{corollary}
 \begin{proof} By Corollary \ref{desc}, the descendant subdigraph of $D$ satisfies P0, P1 and P2. The result then follows from Proposition \ref{conditionstar}. \end{proof} 
 
 \begin{corollary}\label{conditiongamma} Let  $\Gamma$ be a digraph satisfying P0 to P3. Then $\mathcal{C}_{\Gamma}(\alpha)$ holds. In particular, the digraph $\Gamma \setminus \{\alpha\}$,  obtained from $\Gamma$ by deleting the root $\alpha$,  has at least two infinite components.
\end{corollary}
 \begin{proof} By Proposition \ref{conditionstar}, there is $x\in \Gamma$ such that $\mathcal{C}_{\Gamma}(x)$ holds. By Property P1, $\Gamma$ is isomorphic to $\Gamma(x)$. Hence $\mathcal{C}_{\Gamma}(\alpha)$ holds. Let $U,V$ be non-empty subsets of $\Gamma^1$ with ${\rm desc}(U) \cap {\rm desc}(V)=\varnothing$. Then ${\rm desc}(U)$ and ${\rm desc}(V)$ are two distinct infinite components of  $\Gamma \setminus \{\alpha\}$.
  \end{proof}

 
 


 \begin{theorem} \label{disjoint}Let $\Gamma$ be a digraph satisfying  P0 to  P3, and  let $A= {\rm Aut}(\Gamma)$. There is a system $\Omega=\{\omega_1,\ldots, \omega_s\}$ of blocks of imprimitivity for $A$ in  $\Gamma^1$  with $s\geq 2$ and ${\rm desc}(\omega_i) \cap {\rm desc}(\omega_j)=\varnothing$ for $i \neq j$.
\end{theorem}
 \begin{proof} 
 Define a relation on $\Gamma^1$ by saying that vertices $u$ and $v$ are related if, and only if,  $u$ and $v$ belong to the same connected component of the digraph $\Gamma \setminus \{\alpha\}$ obtained from $\Gamma$ by deleting the root $\alpha$.  This is an equivalence relation on $\Gamma^1$ which is invariant under $A$ and, therefore, its equivalence classes form a system of blocks for $A$. Now, by Corollary \ref{conditiongamma}, $\mathcal{C}_{\Gamma}(\alpha)$ holds.  So there are non-empty subsets  $U,V$ of  $\Gamma^1$ with  ${\rm desc}(U) \cap {\rm desc}(V)=\varnothing$.  Since there is no path in $\Gamma \setminus \{\alpha\}$ from a vertex in $U$ to a vertex in $V$,  there are at least two distinct equivalence classes in $\Gamma^1$. The  result follows.
  \end{proof}

  \begin{corollary} \label{disjoint2}Let $D$ be a distance-transitive digraph of  finite out-valency $m>0$ with no directed cycles, and let $G={\rm Aut}(D)$. Suppose the descendant subdigraph $\Gamma={\Gamma}(\alpha)$ of $D$ satisfies  P3. Then the following holds.
  \begin{enumerate}[a)]
\item There is a system  $\Omega=\{\omega_1,\ldots, \omega_s\}$ of blocks of imprimitivity  for $G_{\alpha}$ in $\Gamma^1$  with $s\geq 2$ and ${\rm desc}(\omega_i) \cap {\rm desc}(\omega_j)=\varnothing$ for $i \neq j$.
\item  The permutation group $G^{\Gamma}_{\alpha}$ induced by the action of $G_{\alpha}$ on $\Gamma$ is the wreath product $H\wr_{\Omega}T$, where $H$ is the permutation group $G^{{\rm desc}(\omega_1)}_{\alpha}$   and $T$ is a transitive subgroup of ${\rm Sym}(\Omega)$, the symmetric group on $\Omega$.
\end{enumerate} 
\end{corollary}
  \begin{proof}
  (a) By Corollary \ref{desc},  $\Gamma$ satisfies P0 to P2. The result then follows from Theorem \ref{disjoint}.
  
 (b) Consider the natural action of $G_{\alpha}$ on   $\Omega$. By edge transitivity of $D$, the permutation group $G^{\Omega}_{\alpha}$ is a transitive subgroup $T$ of the symmetric group ${\rm Sym}(\Omega)$. Also, since for $i \neq j$ the sets  ${\rm desc}(\omega_i)$ and  $ {\rm desc}(\omega_j)$ do not intersect and the subdigraph induced on ${\rm desc}(\omega_i)$ is isomorphic to the subdigraph induced on ${\rm desc}(\omega_j)$,  we conclude that the pointwise stabiliser of $\Omega$ in $G_{\alpha}^{\Gamma}$ is the direct product $H  \times \ldots \times H$  of $s$ copies of the permutation group $H:=G^{{\rm desc}(\omega_1)}_{\alpha}$.  The result follows.
   \end{proof}

 \begin{remark}Let $D$ be a primitive, distance-transitive digraph of  finite out-valency $m>0$  with no directed cycles. By Corollary \ref{desc}, $m>1$ and  the descendant subdigraph  $\Gamma$ of $D$ satisfies P0 to P3. Therefore the statement in  Corollary \ref{disjoint2} holds for $D$.
\end{remark}

\subsection{Out-valency equal to a product of two primes}\label{primes}
We now restrict our attention to digraphs satisfying P0 to P3 with  out-valency equal to a product of two primes. For $n\geq 1$, we let $\delta_n$ be the relation as in Definition \ref{delta}, and let $N$ be the number as in Definition \ref{defN}. 
\begin{proposition} \label{pq}Let $\Gamma$ be a digraph  satisfying  P0 to  P3. Suppose the relation $\delta_{N-1}$ is trivial. If the out-valency $m$ is the product of two primes $p$ and $q$, with $p\leq q$,  then the following holds.
 \begin{enumerate}[a)]
\item If $p=q$, then $\Gamma$ is a (rooted) tree.
\item If $p< q$ and $\Gamma$ is not a tree, then there is a system $\{\omega_1,\ldots, \omega_s\}$ of blocks  for $A={\rm Aut}(\Gamma)$ in  $\Gamma^1$  with $s=p$,  $ |\omega_i|=q$ for all $i$ and ${\rm desc}(\omega_i) \cap {\rm desc}(\omega_j)=\varnothing$ for $i \neq j$. Moreover, $r_i=1$ for $i \leq N-1$, and $r_i=p$ for $i \geq N$.
\end{enumerate} 
\end{proposition}
\begin{proof}
By Theorem \ref{disjoint}, there is a system of blocks  $\{\omega_1,\ldots, \omega_s\}$ for  $A$ in $\Gamma^1$ with ${\rm desc}(\omega_i) \cap {\rm desc}(\omega_j)=\varnothing$ for $i \neq j$. Since the cardinality  of $\Gamma^1$ is the product $p\times q$, the number $s$ of blocks is either equal to $p$ or equal to $q$. We now consider the possibilities for the value of the ultimate in-valency $r_N$ of  $\Gamma$. By Lemma \ref{ultimate},  $r_N$ divides the out-valency $m$. Also, by Lemma \ref{ri},  $r_N < m$, since $\Gamma$ satisfies P3. So we have three possibilities: $r_N=1$, $ r_N=p$ or $r_N=q$. Note that if $r_N=1$ then $r_i=1$ for all $i$ and it follows that $\Gamma$ is a tree.

Consider the level $\Gamma^N$. This is the disjoint union of the sets $\Gamma^{N-1}(\omega_1),\ldots, \Gamma^{N-1}(\omega_s)$ where $| \Gamma^{N-1}(\omega_i)|=|\Gamma^{N-1}(\omega_j)|$ for all $i,j$. Then 
\begin{equation}\label{1}
|\Gamma^N|=\sum_{i=1}^s | \Gamma^{N-1}(\omega_i)|=s \times | \Gamma^{N-1}(\omega_1)|.
\end{equation}
To ease notation let $\omega:=\omega_1$. For $x \in  \omega$, $ | \Gamma^{N-1}(x)| \leq | \Gamma^{N-1}(\omega)|$ since $ \Gamma^{N-1}(x)\subseteq  \Gamma^{N-1}(\omega)$. Also, by P1, $ | \Gamma^{N-1}(x)| = | \Gamma^{N-1}|$. This together with (\ref{1}) gives us that 
\begin{equation}\label{2}
|\Gamma^N|=s \times | \Gamma^{N-1}(\omega)| \geq s \times | \Gamma^{N-1}|.
\end{equation}
By replacing $|\Gamma^N|=| \Gamma^{N-1}|  \times \frac{m}{r_N}$ in the  inequality (\ref{2}) and  then dividing by $s$, we obtain 
\begin{equation}\label{3}
| \Gamma^{N-1}|  \times \frac{pq}{s\times r_N}= | \Gamma^{N-1}(\omega)| \geq  | \Gamma^{N-1}|.
\end{equation}

First assume that  $p = q$. In this case we claim that $r_N = 1$ and therefore $\Gamma$ is a tree.  Seeking a contradiction, suppose $r_N \neq 1$. Then the only possibility is that $s=r_N=p$. Replacing these values in (\ref{3}), we obtain $| \Gamma^{N-1}| = | \Gamma^{N-1}(\omega)|$. It then  follows that $\Gamma^{N-1}(x) =  \Gamma^{N-1}(y)$ for all $x,y\in \omega$, since $\Gamma^{N-1}(x) \subseteq  \Gamma^{N-1}(\omega)$ and $| \Gamma^{N-1}(x)| = | \Gamma^{N-1}|$. This contradicts the assumption that the relation  $\delta_{N-1}$ is trivial since the block $\omega$ contains at least 2 vertices. Hence $r_N=1$ and $\Gamma$ is a tree. This proves (a).

Now assume that $p < q$ and that $\Gamma$ is not a tree. Then $r_N \neq 1$. Suppose  $s=q$. Then, by (\ref{3}), $r_N \leq p$. So $r_N=p$ since $r_N$ divides $m$ and  $r_N\neq 1$. Again, by  (\ref{3}), we conclude that $| \Gamma^{N-1}| = | \Gamma^{N-1}(\omega)|$, and this contradicts the assumption that $\delta_{N-1}$ is trivial.  Hence $s \neq q$. We now know that $s=p$ and, by (\ref{3}), $r_N \leq q$. If $r_N=q$, once again we obtain that $\delta_{N-1}$ is not trivial. So the only possibility is that $r_N=p$. It then follows from the definition of $r_N$ and the fact that $r_i$ divides $p^iq^i$, that $r_i=1$ for $i \leq N-1$ and $r_i=p$ for $i \geq N$.
\end{proof}

\vspace{0.3cm}

We are now ready to prove Theorem \ref{prime2}.
 
 \noindent\textbf{Proof of Theorem \ref{prime2}.\/} Let $D$ be a primitive, distance-transitive digraph with no directed cycles and out-valency $p^2$, where $p$ is a prime. By Corollary \ref{desc}, the descendant subdigraph $\Gamma$ of $D$ satisfies P0 to P3. Also, by primitivity of ${\rm Aut}(D)$,  the relation $\delta_{n}$ is trivial for all $n$. The result then follows from Proposition \ref{pq}(a). \hfill $\Box$ 

\begin{corollary} \label{pq2}Let $\Gamma$ be a digraph  satisfying  P0 to  P3, and let $A={\rm Aut}(\Gamma)$. Suppose the relation $\delta_{N-1}$ is trivial. If the out-valency of $\Gamma$ is the product of two primes $p$ and $q$,  with $p\leq q$, then there is a subgroup $H$ of the permutation group $A^{ \Gamma^1}$ such that $H \cong \mathbb{Z}_q \wr_{\Omega} {\rm Sym} (\Omega)$, where  ${\rm Sym} (\Omega)$ is the symmetric group on a set $\Omega$ of cardinality $p$.
\end{corollary}
\begin{proof} If  $\Gamma$ is a tree then the result holds. Suppose $\Gamma$ is not a tree. By Proposition \ref{pq},  $p\neq q$ and there is a system of blocks $\Omega=\{\omega_1,\ldots, \omega_p\}$ for $A$ in  $\Gamma^1$  with blocks of cardinality $q$ and ${\rm desc}(\omega_i) \cap {\rm desc}(\omega_j)=\varnothing$ for $i \neq j$. Since ${\rm desc}(\Gamma^1)$ is the disjoint union of the sets ${\rm desc}(\omega_1), \ldots , {\rm desc}(\omega_p)$, the permutation group $A^{\Omega}$ is the full symmetric group ${\rm Sym}(\Omega)$.

Now  let $i \in \{ 1,2, \ldots, p \}$ and consider the permutation group $A^{ \omega_i}$. Since  $|\omega_i|=q$ is a prime, there is a subgroup of $A^{ \omega_i}$  which is isomorphic to $\mathbb{Z}_q $. There are $p$ such blocks in $\Gamma^1$ and any two distinct blocks have disjoint descendant sets. So the direct product $\mathbb{Z}_q  \times \ldots \times \mathbb{Z}_q$  of $p$ copies of $ \mathbb{Z}_q$ is a subgroup of the pointwise stabiliser of $\Omega$ in $A^{\Gamma^1}$. The result follows.
\end{proof}

 By Proposition \ref{pq} and Corollary \ref{pq2}, we conclude the following.
 \begin{theorem}\label{conclusion} Let $D$ be a distance-transitive  digraph of finite out-valency $m>0$ with no directed cycles, and let $G={\rm Aut}(D)$.  Suppose its descendant subdigraph $\Gamma={\Gamma}(\alpha)$ satisfies P3 and $\delta_{N-1}$ is trivial. If the out-valency of $D$ is a product of two primes $p$ and $q$,  with $p\leq q$,  then there is a system $\Omega=\{\omega_1,\ldots, \omega_p\}$ of blocks for $G_{\alpha}$ in  $\Gamma^1$ with $ |\omega_i|=q$ for all $i$ and ${\rm desc}(\omega_i) \cap {\rm desc}(\omega_j)=\varnothing$ for $i \neq j$. Moreover, the ultimate in-valency $r_N=p$ and there is a subgroup $H$ of the permutation group $G_{\alpha}^{ \Gamma^1}$ such that $H \cong \mathbb{Z}_q \wr_{\Omega}\mathbb{Z}_p$.
 \end{theorem}

 \section{Distance-transitive, weakly descendant-homogeneous digraphs} \label{imprimitive}
 In this section we consider distance-transitive digraphs which are weakly descendant-homogeneous. Recall that a digraph $D$ is weakly descendant-homogeneous} if it is vertex transitive and for any two finite subsets $X$ and $Y$  of vertices and an isomorphism $f$ from  ${\rm desc}(X)$ to  ${\rm desc}(Y)$, there is an automorphism of $D$ agreeing with $f$ on $X$. 
 
 Let $D$ be a distance-transitive, weakly descendant-homogeneous digraph of finite out-valency with no directed cycles. In Theorem \ref{hasP3} we show that if the descendant subdigraph of $D$ has property P3 then $D$ is a digraph of infinite in-valency, $D$ does not have property $Z$ and the reachability relation  $\mathcal{A}$ is universal. We then use this result to prove Theorem \ref{locallyfiniteD}.
 Continuing, in Section \ref{teorema1.3} we show that there is $n \geq 1$ such that $\delta_n$ is non-trivial, where $\delta_n$  is the relation in Definition \ref{delta}.  As there is no $n$ such that $\delta_n$ is universal, otherwise $D$ contains a directed cycle, we conclude that the digraph $D$ is imprimitive. Finally, we investigate the quotient digraph $D/\delta_{n}$ and prove Theorem  \ref{equivalence1}.

 \subsection{Proof of Theorem \ref{locallyfiniteD}} \label{teorema1.2}
 
 
  
  
  The following result appears in Section 3 of  {\cite{AT2}}  and gives a sufficient condition for a weakly descendant-homogeneous digraph not to have property $Z$. 
  Recall that a digraph $D$ has property $Z$ if there is a homomorphism from $D$ onto $Z$, where $Z$ is the digraph which has the set of integers $\mathbb{Z}$ as its vertex set and edges $(i,i+1)$ for every $i \in \mathbb{Z}$. 
  
\begin{lemma}({\cite[Lemma 3.1]{AT2}}) \label{propZ1}
Let $D$ be a weakly descendant-homogeneous digraph, and suppose that for any $x \in D$, there are vertices $a,b \in {\rm desc}(x)$ such that ${\rm desc}(a) \cap {\rm desc}(b)=\varnothing$. Then $D$ is connected and does not have property $Z$.
\end{lemma}

We recall the definition of the reachability relation $\mathcal{A}$. If $e$ and $e'$ are edges in $D$ and there is an alternating walk $x_0x_1\ldots  x_n$ such that $e=(x_0,x_1)$ and either $e'=(x_{n-1},x_n)$ or $e'=(x_n,x_{n-1})$, then $e'$ is said to be {\itshape reachable} from $e$ by an alternating walk and this is denoted by $e \mathcal{A} e'$. Clearly $\mathcal{A}$ is an equivalence relation on the edge set of $D$ which is preserved by the automorphism group of $D$. For an edge $e$ in $D$, we denote its equivalence class by $\mathcal{A}(e)$. 

 \begin{proposition} \label{infinitedegree} Let $D$ be a digraph such that its descendant subdigraph $\Gamma$  satisfies P0 to P2. Suppose there are $x \in \Gamma$ and distinct vertices  $a,b \in \Gamma^1(x)$ such that ${\rm desc}(a) \cap {\rm desc}(b)=\varnothing$.  If $D$ is weakly descendant-homogeneous then $D$ has infinite in-valency, the reachability relation $ \mathcal{A}$ is universal and $D$ does not have property $Z$.
 \end{proposition}
 \begin{proof} 
 By assumption, there is $x \in \Gamma$ and $a,b \in \Gamma^1(x)$ such that ${\rm desc}(a) \cap {\rm desc}(b)=\varnothing$. Consider an infinite arc $bb_1 b_2  \ldots b_i \ldots$  starting at $b$. Then, for $i \geq 1$, ${\rm desc}(a) \cap {\rm desc}(b_i)=\varnothing$ and there is an isomorphism $f_i$ from  the subdigraph induced on ${\rm desc}(a) \cap {\rm desc}(b)$ to the subdigraph induced on ${\rm desc}(a) \cap {\rm desc}(b_i)$ which fixes $a$ and takes $b$ to $b_i$. For $i \geq 1$ there is an automorphism $g_i$ of $D$ which agrees with $f_i$ on $\{ a,b \}$, since $D$ is weakly descendant-homogeneous.  We now show that $x^{g_i} \neq x^{g_j}$ for all $i \neq j$. To ease notation, for  $i\geq 1$ we let $x_i:=x^{g_i}$.
  
 Seeking a contradiction, suppose $x_i= x_j$ for some $i \neq j$. Without loss of generality, assume $j >i$. First note that for $l \geq 1$, the image $b_l=b^{g_l}$ lies in  $\Gamma^1(x_l)$ since $b  \in \Gamma^1(x)$. Then $b_j$ lies in   $\Gamma^{j-i+1}(x_i)$ since $b_j \in \Gamma^{j-i}(b_i)$ and $b_i \in \Gamma^1(x_i)$. Also, $b_j \in \Gamma^1(x_i)$ since $b_j \in \Gamma^1(x_j)$ and we are assuming $x_i=x_j$. So $b_j$ lies in the intersection $\Gamma^{j-i+1}(x_i)\cap  \Gamma^{1}(x_i)$ and this contradicts P0 since $i \neq j$.  Hence $ \{x, x_1,x_2, \ldots  \}$ is an infinite set of vertices which is contained in the set of in-neighbors of $a$. Thus $D$ has infinite in-valency.
 
 To show that the reachability relation is universal, consider the $2$-arc $ x b b_1$. The path $bxax_1 b_1b$
  is an alternating walk from the edge $e=(x,b)$ to the edge $f=(b,b_1)$. Hence $e \mathcal{A} f$ and it follows that the equivalence class $ \mathcal{A}(e)$ contains a $2$-arc. By \cite[Proposition 1.1(a)]{Cameronetal}, the fact that  $ \mathcal{A}(e)$ contains a $2$-arc implies that $ \mathcal{A}$ is universal. 
 
 Finally, by our assumption on $\Gamma$,  there are vertices $a,b$ in $\Gamma$ such that ${\rm desc}(a) \cap {\rm desc}(b)=\varnothing$. So the fact that $D$ does not have property $Z$  follows from Lemma  \ref{propZ1}.
 \end{proof} 
 
 \vspace{0.3cm}
 
From  Propositions \ref{conditionstar}  and \ref{infinitedegree} we obtain the following.
 
  \begin{theorem} \label{hasP3}Let $D$ be a distance-transitive, weakly descendant-homogeneous digraph of finite out-valency $m>0$ with no directed cycles. Assume that its descendant subdigraph satisfies  P3. Then $D$ has infinite in-valency, the reachability relation $ \mathcal{A}$ is universal and $D$ does not have property $Z$. 
 \end{theorem}

 We now prove Theorem \ref{locallyfiniteD}. 
 
\noindent\textbf{Proof of Theorem \ref{locallyfiniteD}.\/} Let $\Gamma$ be the descendant subdigraph of $D$. If $\Gamma$ satisfies P3 then, by Theorem \ref{hasP3},  $D$ has infinite in-valency. This contradicts the assumption that $D$ is locally finite. Hence $\Gamma$ does not satisfy P3 and the result follows from Corollary \ref{finitevalency}. \hfill $\Box$

 \subsection{The quotient digraph} \label{teorema1.3}
 

 Let $n \geq 1$ and let $\delta_n$ be  as  in Definition \ref{delta}. 
 For  a digraph $D$ and $u \in D$, we denote by $\bf{u}$ the $\delta_n$-class containing $u$, and by ${\bf{D}}_n$ the quotient digraph $D/\delta_{n}$. 
  In this section we investigate the quotient digraph ${\bf{D}}_n$ and prove Theorem \ref{equivalence1}.

  \begin{theorem} \label{primitive}Let $D$ be a distance-transitive, weakly descendant-homogeneous digraph of  finite out-valency $m>0$ with no directed cycles. There is $n\geq 1$ such that $\delta_n$ is non-trivial.
  \end{theorem}
   \begin{proof} 
  Let $\Gamma=\Gamma(\alpha)$ be the descendant subdigraph of $D$. By Corollary \ref{desc}, the digraph $\Gamma$ satisfies P0, P1 and P2. Now, by Theorem \ref{regular},  there is an integer $k$ such that for all $l \geq 2k-1$ the subdigraph induced on the descendant set  ${\rm desc}(\Gamma^l(\alpha))$ of the layer $\Gamma^l(\alpha)$ (where $\alpha$ is the root of $\Gamma$) is the disjoint union of finitely many copies of a fixed  infinite digraph $T^*$, where $T^*$  is isomorphic to $T( {\bf{w}})$ and ${\bf{w}}$ is a $\rho$-class in $\cup_{l\geq2k-1}\Gamma^l(\alpha)$. Recall that $T({\bf{w}})$ is the structure consisting of the induced digraph on $\Gamma({\bf{w}})$ together with the equivalence relation induced by $\rho$. 
  We know that $T( {\bf{w}})/\rho$ is a rooted directed tree of out-valency $s \geq1$. If $s=1$ then $T( {\bf{w}})/\rho$ is a directed line, and it follows from Proposition \ref{equivalence} that $\Gamma$ does not satisfy P3. Hence, by Lemma \ref{invalency}(c), we conclude that $\delta_j$ is non-trivial for  some $j \geq N$, where $N$ is as in Definition \ref{defN}.

   
   Now assume $s \geq 2$. In this case, $\Gamma /\rho$  is not a directed line  and it follows from Proposition \ref{equivalence} that $\Gamma$ satisfies P3. (The remainder of the argument is as in the proof of  \cite[Theorem 3.3]{AT2}). Let $l \geq 2k-1$ and let $T_1^{\ast}, \ldots ,T_n^{\ast}$ be the distinct copies of $T^{\ast}$ in ${\rm desc}(\Gamma^l(\alpha))$. For $i \in \{1, \ldots ,n\}$, let $\varphi_i$ be a homomorphism from $T_i^{\ast}$ onto the infinite rooted regular tree $T$ of out-valency $s$. Note that this homomorphism exists since, as mentioned above, $T^*$  is isomorphic to $T( {\bf{w}})$ and $T( {\bf{w}})/\rho$ is a rooted directed tree of out-valency $s$. 
Let $r\geq2$ and let $\{t_1, \ldots ,t_{s^r}\}$ be the set $T^r(a)$ of vertices at distance
$r$ from the root $a$ of $T$. So the layer $\Gamma^{l+r}(\alpha)$ is the disjoint union 
\[
U:=\bigcup_{i=1}^{n}\left(
\bigcup_{j=1}^{s^r}\varphi_i^{-1}(t_j)\right) . 
\]
Now let $x, y, z \in V(T)$ be such that $x$ and $y$ have a common predecessor $u$ in $T$ but $z$ does not have a
common predecessor with $x$ (and hence not with $y$ either). So $\varphi_1^{-1}(x)$ and $\varphi_1^{-1}(y)$
have a common predecessor $\varphi_1^{-1}(u)$ in $T_1^{\ast}$ but $\varphi_1^{-1}(z)$ does not have a common predecessor with
$\varphi_1^{-1}(x)$ (and hence not with $y$ either). Let $S = \varphi_1^{-1}(y) \cup \varphi_1^{-1}(z)$. We
know that the subdigraph $\Delta$ of $\Gamma$ on the subset desc$(U)$ is the union of finitely many
disjoint copies of $T^{\ast}$, so there is an automorphism $\theta$ of $\Delta$ which fixes
desc$(U\backslash S)$ pointwise and interchanges desc$ (\varphi_1^{-1}(y))$ and desc$(\varphi_1^{-1}(z))$. By weak
descendant-homogeneity, $\theta$ extends to an automorphism $g$ of $D$ which agrees with $\theta$ on $\Gamma^{n+r}(\alpha)$. Since $\varphi_1^{-1}(x)$ and
$\varphi_1^{-1}(y)$ have a common predecessor in $T_1^{\ast}$ and $g$ fixes $T_1^{\ast}$ pointwise,
$\varphi_1^{-1}(x)$ and $\varphi_1^{-1}(y)^{g}(=\varphi_1^{-1}(z)$) have a common predecessor, namely
$\varphi_1^{-1}(u)^g$. Note that this is  the image of the common predecessor, $\varphi_1^{-1}(u)$, of
$\varphi_1^{-1}(x)$ and $\varphi_1^{-1}(y)$ in $T_1^{\ast}$. But $\varphi_1^{-1}(u)$ is not contained in
$\Gamma$ since $\varphi_1^{-1}(x)$ and $\varphi_1^{-1}(z)$ do not have a common predecessor within $\Gamma$.
So $\alpha^{g}\neq\alpha$. Now, as $g$ leaves the set $\Gamma^{l+r}(\alpha)$ invariant, it
follows that $\Gamma^{l+r}(\alpha ) = \Gamma^{l+r}(\alpha^{g})$. Hence, the equivalence relation $\delta_{l+r}$ is non-trivial in the digraph $D$.
    \end{proof}
    
      \begin{corollary} \label{prim}A  distance-transitive, weakly descendant-homogeneous digraph of  finite out-valency $m>0$ with no directed cycles is imprimitive. 
  \end{corollary}
  \begin{proof}
  By Theorem \ref{primitive}, there is $n \geq 1$ such that the equivalence relation $\delta_n$ is non-trivial. If $\delta_n$ is universal then $D$ contains a finite directed cycle, and this contradicts our assumption on $D$. Hence $\delta_n$ is a proper non-trivial equivalence relation which is preserved by the automorphism of $D$. Thus $D$ is imprimitive.
  \end{proof}

  \begin{lemma} \label{transitive}Let $D$ be a distance-transitive, weakly descendant-homogeneous digraph of finite out-valency $m>0$ with no directed cycles. For $n\geq 1$ the quotient digraph  ${\bf{D}}_n$  is infinite and edge transitive.  
 \end{lemma}
\begin{proof} 
 First we show that there are infinitely many distinct $\delta_{n}$-classes. Let $L= \dots v_{-2}v_{-1}v_0\dots v_i\dots$ be a directed line in $D$. We claim  that for $i \neq j$, the classes ${\bf{v}}_i$ and ${\bf{v}}_j$ are distinct. Without loss of generality, suppose $j > i$.  Then $v_j \in {\rm desc}^{j-i}(v_i)$, and it follows that ${\rm desc}^n(v_j)\subseteq {\rm desc}^{n+(j-i)}(v_i)$. If ${\bf{v}}_i={\bf{v}}_j$ then  ${\rm desc}^n(v_i)={\rm desc}^n(v_j)$, and it follows that  ${\rm desc}^n(v_i)\subseteq {\rm desc}^{n+(j-i)}(v_i)$. This contradicts P0 since $j-i >0$. 
   
The fact that ${\bf{D}}_n$ is edge transitive follows from edge transitivity of $D$. Indeed, let $(\bf{a},\bf{b})$ and $(\bf{c},\bf{d})$ be edges in ${\bf{D}}_n$.  By definition of adjacency, there are vertices $x \in \bf{b}$ and $y \in \bf{d}$ such that $(a,x)$ and $(c,y)$ are edges in $D$. By edge transitivity of $D$, there is an automorphism $g$ of $D$ such that $a^g=c$ and $x^g=y$. Since the set of $\delta_{n}$-classes is a system of blocks of imprimitivity for the automorphism group of $D$, the element $g$ takes  $\bf{a}$ to $ \bf{c}$ and takes $\bf{b}$ to $ \bf{d}$. 
 \end{proof}

 \begin{lemma}\label{propZ} Let $D$ be an edge transitive digraph and let  $n\geq 1$. If ${\bf{D}}_n$ has property $Z$ then $D$ has property $Z$.    
  \end{lemma}   
    \begin{proof} 
   Suppose ${\bf{D}}_n$ has property $Z$ and  let  $\phi: {\bf{D}}_n \rightarrow Z$ be a surjective digraph homomorphism. Now consider the projection $\pi: D \rightarrow {\bf{D}}_n$ which takes a vertex $a$ to its class $\bf{a}$. Then the composition $\phi \circ \pi$ is a digraph homomorphism from $D$ onto $Z$. Hence, $D$ has property $Z$
     \end{proof}

 We consider the case when the descendant subdigraph $\Gamma$ does not satisfy P3. 
    
  \begin{lemma} \label{Dtree}Let $D$ be a distance-transitive, weakly descendant-homogeneous digraph of finite out-valency $m>0$ with  no directed cycles. Suppose the descendant subdigraph $\Gamma$ does not satisfy  P3, and let $N$ be as in Definition \ref{defN}. Then, for  $n \geq {N-1}$, the quotient digraph ${\bf{D}}_{n}$ has out-valency 1. 
  \end{lemma}
  \begin{proof} Let $n \geq {N-1}$ and let ${\bf{u}} \in {\bf{D}}_{n}$. Seeking a contradiction, suppose there are distinct vertices $ \bf{v}_1,\bf{v}_2$ in $ {\bf{D}}_{n}$ such that  $(\bf{u},\bf{v}_1)$ and $(\bf{u},\bf{v}_2)$ are edges. By definition of adjacency, there are vertices $a,b \in \bf{u}$ such that $(a,v_1)$ and $(b,v_2)$ are edges in $D$. Since $ v_1\in {\rm desc}^{1}(a)$, we apply  Lemma \ref{invalency}(c) with $t=1$ and $l=n$ to conclude that 
  \begin{equation}\label{n}
  {\rm desc}^n(v_1)={\rm desc}^{n+1}(a).
  \end{equation}
  Similarly, 
  \begin{equation}\label{n+1}
  {\rm desc}^n(v_2)={\rm desc}^{n+1}(b).
    \end{equation}
    This is illustrated in Figure \ref{T1}.
    
    \begin{figure}[htbp]
\begin{center}
    \setlength{\unitlength}{4mm}
\linethickness{0.25mm}
\begin{center}
\begin{picture}(12,10)(-1,1)

\put(5.5,10){\oval(6,1.8)}
\put(2,10.5){\makebox(0,0)[t]{$\textbf{u}$}}
\put(4,10){\circle*{0.4}}
\put(4,10){\vector(-1,-1){1.5}}

\put(7,10){\circle*{0.4}}
\put(7,10){\vector(1,-1){1.5}}

\qbezier(2.5,8.5)(4.5,4.5)(4.5,2)

\put(2.5,8.5){\circle*{0.3}}
\put(8.5,8.5){\circle*{0.3}}

\qbezier(2.5,8.5)(4.5,4.5)(4.5,2)
\qbezier(8.5,8.5)(10.5,4.5)(10.5,2)
\qbezier(0.5,2)(2.5,2)(4.5,2)
\qbezier(8.8,2)(3.5,2)(10.5,2)

\qbezier(2.5,8.5)(0.5,4.5)(0.5,2)
\qbezier(8.5,8.5)(6.5,4.5)(6.5,2)

\put(4,10.8){\makebox(0,0)[t]{$a$}}
\put(7,11){\makebox(0,0)[t]{$b$}}
\put(2,9.35){\makebox(0,0)[t]{$v_1$}}
\put(2.5,1.35){\makebox(0,0)[t]{${\rm desc}^n(v_1)$}}
\put(8.5,1.35){\makebox(0,0)[t]{${\rm desc}^n(v_2)$}}
\put(9,9.35){\makebox(0,0)[t]{$v_2$}}
\end{picture}
\end{center}
\caption{The classes $\textbf{u},\textbf{v}_1$ and $\textbf{v}_2$.}
\label{T1}
\end{center}
\end{figure}
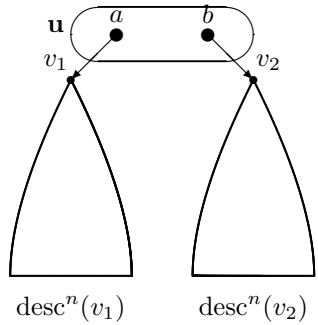

  On the other hand, we have $a,b \in \bf{u}$, and this means that $\delta_n(a,b)$ holds. So ${\rm desc}^n(a)={\rm desc}^{n}(b)$ and it follows that   ${\rm desc}^{i}(a)={\rm desc}^{i}(b)$ for all $i \geq n$. In particular,
  \begin{equation}\label{n+2}
  {\rm desc}^{n+1}(a)={\rm desc}^{n+1}(b).
   \end{equation}
 From (\ref{n}),  (\ref{n+1}) and (\ref{n+2}), we conclude that 
 $$
 {\rm desc}^n(v_1)={\rm desc}^{n}(v_2),
 $$ 
 and this means that  $\delta_n(v_1,v_2)$ holds. This contradicts the assumption that $ \bf{v}_1$ and $\bf{v}_2$ are distinct vertices in  ${\bf{D}}_{n}$. 
 Thus ${\bf{D}}_{n}$ has out-valency 1.
 \end{proof}
   
  \vspace{0.5cm} 
By Lemmas \ref{transitive} and \ref{Dtree},  for  $n  \geq {N-1}$  the quotient digraph   ${\bf{D}}_{n}$ is an infinite  regular directed tree whenever the descendant subdigraph of $D$ does not satisfy P3. 
This fact, together with the next lemma,  will allow us to prove Theorem  \ref{equivalence1}. 
 
 \begin{lemma}\label{cara} Let $D$ be a distance-transitive, weakly descendant-homogeneous digraph of finite out-valency $m>0$ with  no directed cycles. Then $D$ has property $Z$ if, and only if, the descendant subdigraph $\Gamma$ does not satisfy  P3. 
  \end{lemma}
  \begin{proof}Suppose $D$ has property $Z$. Seeking a contradiction, assume $\Gamma$ satisfies P3. By Proposition \ref{conditionstar}, there are vertices $a,b$ in $\Gamma$ such that ${\rm desc}(a) \cap {\rm desc}(b)=\varnothing$. It then follows from Lemma \ref{propZ1} that $D$ does not have property $Z$. 
  
    Conversely, suppose $\Gamma$ does not satisfy P3 and let $n  \geq {N-1}$. By Lemma \ref{Dtree}, the quotient digraph  ${\bf{D}}_{n}$ is an infinite directed tree and therefore has property $Z$. It then follows from Lemma \ref{propZ} that  $D$ has property $Z$. 
   \end{proof}
 
 \vspace{0.3cm} 
   
\noindent\textbf{Proof of Theorem \ref{equivalence1}.\/}  By Lemma \ref{cara}, (a) and (c) are equivalent. We now show that (a) and (b) are equivalent.
We already know that  if $D$ has property $Z$ then the reachability relation is not universal.
 Conversely,  assume that the reachability relation is not universal. Seeking a contradiction, assume $D$ does not have property $Z$. Then, by Lemma \ref{cara}, the descendant subdigraph $\Gamma$ satisfies P3. Hence, by  Theorem \ref{hasP3}, the reachability relation is universal. \hfill $\Box$

 \vspace{0.3cm} 
 \begin{corollary}
 Let $D$ be a distance-transitive, weakly descendant-homogeneous digraph of finite out-valency $m>0$ with no directed cycles. If $D$ is locally finite then its descendant subdigraph does not satisfy  P3.
 \end{corollary}
  \begin{proof}
  Since $D$ is locally finite, we conclude from Theorem \ref{locallyfiniteD} that $D$ has property $Z$. The result then follows from  Lemma \ref{cara}.
  \end{proof}

  
 \vspace{0.3cm} 
 Finally, we end this section considering the case when the descendant subdigraph $\Gamma$ of $D$ satisfies P3. 
  
\begin{lemma} \label{tree}Let $D$ be a distance-transitive, weakly descendant-homogeneous digraph of finite out-valency $m>0$  and with no directed cycles. Suppose the descendant subdigraph $\Gamma$ satisfies  P3. Then $m>1$ and,  for $n \geq 1$, the quotient digraph ${\bf{D}}_n$ has the following properties: 
 \begin{enumerate}[a)]
\item   ${\bf{D}}_{n}$  is edge transitive and does not have property $Z$;
\item   ${\bf{D}}_{n}$ has infinite in-valency;
\item the out-valency of  ${\bf{D}}_{n}$ is at least $s\geq 2$, where $s$ is the number of blocks in $\Gamma^1$ as in Theorem \ref{disjoint}. In particular, for $n=1$,  ${\bf{D}}_{n}$  has out-valency at most $m$. 
\end{enumerate} 
\end{lemma}
\begin{proof} It is clear that $m>1$, otherwise $\Gamma$ does not satisfy   P3. 
  Let $n \geq 1$. If  $\delta_{n}$ is trivial then   ${\bf{D}}_{n}=D$ and the result holds trivially. Assume $\delta_{n}$ is non-trivial.  
  
(a)  The fact that  ${\bf{D}}_{n}$ is edge transitive was proved in Lemma \ref{transitive}. We now show that ${\bf{D}}_{n}$  does not have property $Z$.  Seeking a contradiction, assume  ${\bf{D}}_{n}$ has property $Z$. By Lemma \ref{propZ},  $D$ has property $Z$. Then, by Lemma \ref{cara},  the descendant subdigraph $\Gamma$ does not satisfy P3. 
  
(b) Since $\Gamma$ is a digraph satisfying P0 to P3, by Proposition \ref{conditionstar} there is a vertex $x \in \Gamma$ such that $\mathcal{C}_{\Gamma}(x)$ holds. This means that there are vertices $a,b \in \Gamma^1(x)$ such that ${\rm desc}(a) \cap {\rm desc}(b)= \varnothing$. Let $bb_1b_2 \ldots b_i \ldots$ be an infinite arc starting at $b$. We have shown in the proof of Proposition \ref{infinitedegree} that there is an infinite set $\{ x, x_1, x_2, \ldots, x_i, \ldots \}$ of vertices such that $(x_i,a), (x_i,b_i)$ are edges in $D$ for all $i$. We claim that for $i \neq j$ the vertices $x_i$ and $x_j$ are not in the same $\delta_{n}$-class. Seeking a contradiction, suppose  there are $i,j$ with $i \neq j$ such that   ${\bf{x}}_i = {\bf{x}}_j$.  This means that $\Gamma^{n}(x_i)=\Gamma^{n}(x_j)$. Without loss of generality, we may assume that $j >i$. We have $b_j  \in \Gamma^{j-i}(b_i)$ and $b_i  \in \Gamma^{1}(x_i)$. So $b_j  \in \Gamma^{1+j-i}(x_i)$, and it follows that the set 
 $\Gamma^{n-1}(b_j)$ is contained in $\Gamma^{n+j-i}(x_i)$. On the other hand, $\Gamma^{n-1}(b_j) \subseteq \Gamma^{n}(x_i)$ since $b_j  \in \Gamma^{1}(x_j)$ and $\Gamma^{n}(x_i)=\Gamma^{n}(x_j)$. Hence $\Gamma^{n-1}(b_j)$ is contained in the intersection $\Gamma^{n+j-i}(x_i) \cap \Gamma^{n}(x_i)$, where $n+j-i>n$. This contradicts P0. Thus $\{{\bf{x}}_1,{\bf{x}}_2,\ldots,{\bf{x}}_i,\ldots \}$  is an infinite set of vertices in ${\bf{D}}_{n}$ such that $({\bf{x}}_i,\bf{a})$ is an edge. 
  
(c) Let $a \in D$, $\Gamma={\rm desc}(a)$ and $G={\rm Aut}(D)$. By Theorem \ref{disjoint}, there is a system of blocks $\{\omega_1, \ldots, \omega_s\}$  for $G_a$ in $\Gamma^1(a)$ with  $s \geq 2$ and ${\rm desc}(\omega_i)\cap {\rm desc}( \omega_j)= \varnothing$ for all $i \neq j$. So for  $x \in \omega_i$ and $y \in \omega_j$, the classes $\bf{x}$ and $\bf{y}$ are distinct whenever $i \neq j$. This shows that there are at least $s$ distinct classes  ${\bf{x}}_1,{\bf{x}}_2,\ldots,{\bf{x}}_s$ with  $x_i \in \omega_i$ such that $({\bf{a}},{\bf{x}}_i)$ is an edge in ${\bf{D}}_{n}$  for $i \in \{1,2, \ldots, s\}$. 

Finally,  we consider the case $n=1$. Seeking a contradiction, assume ${\bf{D}}_{1}$ has out-valency greater than $m$. Then there is a vertex $y \notin \Gamma^1(a)$ with $\bf{y} \neq \bf{x}$ for all $x \in  \Gamma^1(a)$ and such that $({\bf{a}},{\bf{y}})$ is an edge in ${\bf{D}}_{1}$. By definition of adjacency, there are vertices $u \in {\bf{a}}$ and $v\in {\bf{y}}$ such that $(u,v)$ is an edge in $D$. Since $u\in {\bf{a}}$, we have $\Gamma^1(u)=\Gamma^1(a)$, and it follows that $v$ lies in $\Gamma^1(a)$. This  implies that $\bf{v} = \bf{x}$ for some $x \in  \Gamma^1(a)$ and therefore $\bf{y} = \bf{x}$. 
\end{proof}


\section{Classification of some distance-transitive weakly descendant-homogeneous digraphs}  \label{classification}
In this section we investigate further   the class of distance-transitive weakly descendant-homogeneous  digraphs for which the reachability relation is not universal. 
We are particularly interested in the digraphs for which the relation $\delta_1$ is non-trivial. 
\vspace{0.2cm}

A bipartite digraph is a directed graph whose vertices are divided into two disjoint and independent sets $X$ and $Y$ such that all edges are directed from $X$ to $Y$. The sets $X$ and $Y$ are called the parts of the digraph.  Let $\Sigma$ be a bipartite digraph with bipartition $X \cup Y$ and no isolated vertices (that is, vertices with degree zero). Assume further that $\Sigma$ is edge transitive. So it is transitive on both parts $X$ and $Y$, and any two vertices in $X$ have the same out-valency $m>0$ and any two vertices in $Y$ have the same in-valency $r>0$. We will also refer to $X$ as the set of sources of $\Sigma$, and refer to $Y$ as its set of sinks. 

We consider the relation $\delta_1$   on $X$: for vertices $u,v$ in $X$, write $\delta_1(u,v)$  if, and only if, $u$ and $v$ have the same set of out-neighbours. This  is an equivalence relation on $X$ which is invariant under the action of  ${\rm Aut}(\Sigma)$.  Note that $\delta_1$ is universal if, and only if, $\Sigma$ is complete bipartite. To ease notation, throughout we let $\delta:=\delta_1$
 

\vspace{0.2cm}


Let $D$ be a connected, edge transitive digraph. Recall that  the subdigraph spanned by one of the equivalence classes of the reachability relation $\mathcal{A}$  is called an  {\itshape{alternet}}. By edge transitivity,  all the alternets of $D$ are isomorphic to some fixed digraph $\Delta(D)$. If the reachability relation has more than one class then it follows from Proposition 1.1 of {\cite{Cameronetal}}  that the digraph $\Delta(D)$ is bipartite. 
 The {\itshape{digraph of alternets}} of $D$, denoted by ${\rm Al}(D)$, is the digraph that has the set of alternets in $D$ as the vertex set, and if $A$ and $B$ are alternets then $(A,B)$ is an edge in ${\rm Al}(D)$ if the intersection of the set of sinks of $A$ with the set of sources of $B$ is non-empty. In particular, if such intersection is trivial, we say that ${\rm Al}(D)$ is {\itshape{loosely attached}}. 
\begin{definition}\label{Rclass}
Let $D$ be a connected, edge transitive digraph. For vertices $u,v$ in $D$, write   $R(u,v) $ if, and only if, $u$ and $v$  belong to the set of sinks of the same alternet and they also belong to the set of sources of some common alternet. Otherwise, say that $R(u,v)$ does not hold.
\end{definition}

 
 \begin{lemma}  \label{alternets} Let $D$ be an edge transitive, weakly descendant-homogeneous digraph of finite out-valency $m>1$ with property $Z$.  Suppose the relation $\delta$ is non-trivial in $D$. Then for any edge $(A,B)$ in ${\rm Al}(D)$, the intersection of the set of sinks of $A$ with the set of sources of $B$ is the union of $\delta$-classes. 
\end{lemma}
\begin{proof}  
Let $\phi: D \rightarrow Z$ be a surjective digraph homomorphism and, for $i \in \mathbb{Z}$, set $\Gamma_i=\phi^{-1}(i)$, the fibres of $D$.  Let $G={\rm Aut}(D)$. Following a similar reasoning as in Corollary 16 of \cite{Moller5}, consider the relation $R$ in Definition \ref{Rclass} on the vertices of $\Gamma_0$. So an $R$-class is the intersection of the set of sinks of an alternet with the set of sources of some other alternet. This is an equivalence relation which is preserved by the subgroup $G_{\{\Gamma_0\}}$. So, both the set of sinks and the set of sources of any alternet is a union of $R$-classes. In particular,  if the alternets are finite, then the cardinality of the $R$-classes divides the cardinality of the set of sources, and also divides the cardinality of the set of sinks.

For an alternet $A$ in $D$, denote by $X_A$ the set of sources of $A$ and by $Y_A$ its set of sinks. Let $A$ and $B$ be alternets in $D$ such that $(A,B)$ is an edge in ${\rm Al}(D)$. Let  $U$ denote the $R$-class $ Y_A \cap X_B$. We know  that  $X_B$ is a union of $\delta$-classes. Seeking a contradiction, suppose that $U$ is not a union of $\delta$-classes. Then there is a $\delta$-class ${\bf b}$  in $X_B$ such that ${\bf b}\cap U \neq \varnothing $ but  ${\bf b}$ is not completely contained in $U$. Let  $v$ be a vertex in ${\bf b}\setminus U$. Then $v$ lies in the set $Y_C$ of sinks of some alternet  $C \neq A$ with $(C,B)$ an edge in ${\rm Al}(D)$. 
Now let $b\in {\bf b}\cap U$. The set $Y_A$ contains at least two distinct vertices since it  contains the set of out-neighbors of any vertex in $X_A$ and we are assuming $m >1$. So there is $a \in Y_A$ with $a \neq b$. Consider the (finitely generated) subdigraph $W$ induced on ${\rm desc}(\{a,b,v\})$ and let $S=\{a,b,v\} \cup {\rm desc}^1(\{a,b\})$. Since ${\rm desc}^1(b)={\rm desc}^1(v)$, there is an automorphism $f$ of $W$ which interchanges $b$ and $v$, and fixes $S \setminus \{b,v\}$. By weak descendant-homogeneity, there is an automorphism $g$ of $D$ agreeing with $f$ on $S$. So $g$ interchanges  the alternets $A$ and $C$ since $b \in A$ and $v \in C$. On the other hand, $g$ fixes $A$ setwise since $a^g=a$. This is impossible since $A$ and $C$ have no common vertices. Thus ${\bf b}$ must be contained in $U$. 
\end{proof}

 \begin{lemma}  \label{sizes} Let $D$ be an edge transitive, weakly descendant-homogeneous digraph of finite out-valency $m>1$ with property $Z$.  Let $A$ be an alternet in $D$ with bipartition $X \cup Y$. Then $X$ is finite if, and only if, $Y$ is finite. Furthermore, if $D$ is locally finite and distance-transitive, then $ |X| \geq  |Y|$ whenever $X$ is finite.  \end{lemma}
 \begin{proof} 
 Suppose $X$ is finite. Then the number of edges in $A$ is equal to $ |X| \times m$, and this is finite. Since any edge $(u,v)$ in $A$ has $v \in Y$, it follows that $Y$ is finite. 
 
 Now suppose $Y$ is finite and consider the relation $\delta$. The set $X$ is a union of $\delta$-classes and each $\delta$-class ${\bf x}$ corresponds to a unique $m$-subset in $Y$, namely the set of out-neighbors of $x$. Since $Y$ is finite, there are only finitely many $m$-subsets in $Y$, and it follows that there are only finitely many $\delta$-classes in $X$. It is now sufficient to show that the $\delta$-classes are finite. This is clear if $\delta$ is trivial. Assume $\delta$ is non-trivial. By Lemma \ref{alternets}, the set $Y$ is a union of $\delta$-classes and, since $Y$ is finite, such classes must be finite. 
 
 Finally, assume $D$ is locally finite and distance-transitive. Let $r$ denote the in-valency of $D$. Since we are assuming that $D$ has property $Z$, it follows from Theorem \ref{equivalence} that the descendant subdigraph $\Gamma$ of $D$ does not satisfy P3. Then, by Corollary \ref{cor41},  we have $r\geq m$. By counting the edges of $A$ in two different ways, we have $ |X| \times m=|Y| \times  r$. The result follows.  
 \end{proof}

\begin{corollary}   Let $D$ be an edge transitive, weakly descendant-homogeneous digraph of finite out-valency $m>1$ with property $Z$.  Suppose the relation $\delta$ is non-trivial in $D$. Then  the digraph ${\rm Al}(D)$ is not loosely attached.  
\end{corollary}
\begin{proof}  Let $(A,B)$ be an edge in ${\rm Al}(D)$. To prove that ${\rm Al}(D)$ is not loosely attached, we must  show that the intersection of the set of sinks of $A$ with the set of sources of $B$ contains at least two vertices. By assumption, a $\delta$-class contains at least two vertices.  Now, by Lemma \ref{alternets}, the intersection of the set of sinks of $A$ with the set of sources of $B$ is the union of $\delta$-classes. The result follows. 
\end{proof} 

\vspace{0.3cm}

 By Lemma \ref{sizes}, if $D$ is an edge transitive, weakly descendant-homogeneous digraph and  $A$ is an alternet in $D$ with bipartition $X\cup Y$, then either both parts, $X$ and $Y$, are finite or both are infinite. Here we consider only the finite case. 

\begin{definition} We let $\mathcal{C}$ be  the class of bipartite digraphs $\Sigma$ on $X\cup Y$ with no isolated vertices and  satisfying the following properties: 
 \begin{itemize}
 \item $X$ and $Y$ are finite non-empty sets.
 \item $\Sigma$ is edge transitive.
 \item $\delta$ is non-trivial on $X$.
 \item The cardinality of $Y$  is equal to the cardinality of the $\delta$-classes in $X$.
 \end{itemize}
 \end{definition}


For $\Sigma \in \cal{C}$, as noted earlier, edge transitivity implies that  $\Sigma$ is transitive on both parts $X$ and $Y$. So any two vertices in $X$ have the same out-valency $m>0$, and any two vertices in $Y$ have the same in-valency $r>0$. Then, as the cardinality of $Y$  is equal to the cardinality of the $\delta$-classes in $X$, by counting edges from  $X$ to $Y$, we conclude that the in-valency $r$ of $Y$ is equal to $k \times m$, where $k$ is the number of $\delta$-classes in $X$.

\vspace{0.3cm}

We now present a family of highly-arc-transitive, descendant-homogeneous digraphs with property Z whose alternets  belong to the class  $\mathcal{C}$. This family of digraphs was originally constructed  in Section 4.3.3 of {\cite{AT2}} but  the notation $D(m,M)$ was introduced only later in Example 4.4 of \cite{Amato2}.

\begin{example} ({\cite[Example 4.4]{Amato2}})\label{example1}
{\rm Let $m$ be an integer greater than or equal to 1 and let $M$ be either an integer with $m \leq M$, or let $M=\aleph_0$.   Let  $\Omega=\{1,2, \ldots\}$ be a set of cardinality $M$. We construct an infinite, countable digraph $D=D(m,M)$  with out-valency $m$  as follows.

Start with an infinite directed regular tree $T$ of out-valency 1 and in-valency $\binom{M}{m}$ if $M<\aleph_0$, and in-valency $\aleph_0$ if $M=\aleph_0$. Let $D$ have vertex set $T \times \Omega$. In other words, the vertex set of $D$ is obtained from the vertex set of $T$ by replacing each $t \in T$ by a copy $\{t\} \times \Omega$ of $\Omega$. Now, for each $t \in T$, let the $m$-element subsets of $\Omega$ be in one-to-one correspondence with the set ${\rm in}_T(t)$ of in-neighbours of $t$ in $T$. We define adjacency in $D$ by: if $(t_i, t) \in E(T)$ and $t_i$ corresponds to the $m$-element subset $\{a_1, a_2, \ldots, a_m\}$ then $((t_i,c),(t,a_l)) \in E(D)$, for all $l \in \{1, 2,\ldots, m\}$ and all $c \in \Omega$. The digraph $D(m,M)$ is infinite, highly-arc-transitive, descendant-homogeneous (see  Theorem 4.6 of {\cite{AT2}}) and has property $Z$. In particular, the digraph $D(1,M)$ is an infinite regular tree of out-valency 1 and in-valency $M$.

\begin{figure}[htbp]
\begin{center}

\setlength{\unitlength}{4mm}
\linethickness{0.2mm}
\begin{center}
\begin{picture}(12,10)(-4.5,0)
\put(-1,0){\oval(6,1.5)}
\put(-0.5,0.00){\circle*{0.3}}
\put(0.5,0.00){\circle*{0.3}}
\put(-2,0.00){\circle*{0.3}}
\put(-3.00,0.00){\circle*{0.3}}
\put(-1.2,0){\makebox(0,0)[t]{{\small {$\ldots$}}}}
\put(-1,4){\oval(6,1.5)}
\put(-3.00,4.00){\circle*{0.3}}
\put(-2.00,4.00){\circle*{0.3}}
\put(-1.2,4){\makebox(0,0)[t]{{\small {$\ldots$}}}}
\put(-0.5,4.00){\circle*{0.3}}
\put(0.5,4.00){\circle*{0.3}}
\put(-3.00,4.00){\vector(0,-1){4}}
\put(-3,4.00){\line(1,-4){1}}
\put(-2,4.00){\vector(0,-1){4}}
\put(-2,4.00){\line(-1,-4){1}}
\qbezier(-0.5,4)(-1,3.2)(-3,0)
\qbezier(-0.5,4)(-1,2,66)(-2,0)
\qbezier(0.5,4)(0,3.42)(-3,0)
\qbezier(0.5,4)(0,3.2)(-2,0)
\put(2.8,4){\makebox(0,0)[t]{{\small {$\ldots$}}}}
\put(6.5,4){\oval(6,1.5)}
\put(4.5,4.00){\circle*{0.3}}
\put(5.5,4.00){\circle*{0.3}}
\put(7.5,4){\makebox(0,0)[t]{{\small {$\ldots$}}}}
\put(8.5,4.00){\circle*{0.3}}
\put(4.5,4.00){\line(-5,-4){5.1}}
\put(4.5,4.00){\vector(-1,-1){4}}
\put(5.5,4.00){\vector(-3,-2){6}}
\put(5.5,4.00){\line(-5,-4){5}}
\qbezier(8.5,4)(4,2)(-0.5,0)
\qbezier(8.5,4)(3,1.25)(0.5,0)
\put(10.5,4){\makebox(0,0)[t]{{\small {$\ldots$}}}}
\put(-1,8){\oval(6,1.5)}
\put(-3.00,8){\circle*{0.3}}
\put(-2.00,8){\circle*{0.3}}
\put(0,8){\makebox(0,0)[t]{{\small {$\ldots$}}}}

\put(11,8){\makebox(0,0)[t]{{\small {$\ldots$}}}}

\put(1.00,8){\circle*{0.3}}
\put(-3.00,8){\vector(0,-1){4}}
\put(-3,8){\line(1,-4){1}}
\put(-2,8){\vector(0,-1){4}}
\put(-2,8){\line(-1,-4){1}}
\put(1.00,8){\line(-1,-1){4}}
\put(1.00,8){\line(-3,-4){3}}
\put(2.8,8){\makebox(0,0)[t]{{\small {$\ldots$}}}}
\put(6.5,8){\oval(6,1.5)}
\put(4.5,8){\circle*{0.3}}
\put(5.5,8){\circle*{0.3}}
\put(7.5,8){\makebox(0,0)[t]{{\small {$\ldots$}}}}
\put(8.5,8){\circle*{0.3}}
\put(4.5,8){\line(-5,-4){5.1}}
\put(4.5,8){\vector(-1,-1){4}}
\put(5.5,8){\vector(-3,-2){6}}
\put(5.5,8){\line(-5,-4){5}}
\qbezier(8.5,8)(4,6)(-0.5,4)
\qbezier(8.5,8)(3,5.25)(0.5,4)
\put(10.5,4){\makebox(0,0)[t]{{\small {$\ldots$}}}}
\put(4,0){\makebox(0,0)[t]{{\small {$\ldots$}}}}
\end{picture}
\end{center}
\caption{Partial view of the Digraph  $D(2, \aleph_0$).}
\label{T2}
\end{center}
\end{figure}
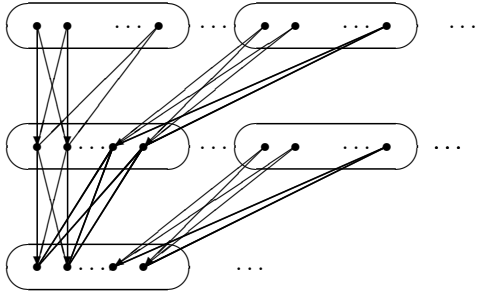

There is clearly a natural homomorphism from $D$ onto $T$ taking $(t,c)$ to $t$ for all $c \in \Omega$. Let $t \in T$. All vertices in $\{t\} \times \Omega$ have the same out-neighbours, namely the set $\{(t',a_1), \ldots, (t',a_m)\}$, where $(t, t')\in E(T)$ and $\{a_1, \ldots, a_m\}$ is the (unique) $m$-element subset of vertices from $\Omega$ corresponding to $t$. So the digraph induced on 
$$
(\{t\}\times \Omega) \cup \{(t',a_1),\ldots, (t',a_m)\}
$$
is the (directed) complete bipartite digraph $\overrightarrow{K}_{M,m}$. Also, the digraph induced on 
$$
(\bigcup_{t_i \in {\rm in}_T(t)}\{t_i\}\times \Omega) \cup (\{t\}\times \Omega)
$$
is an alternet of $D$: this is a bipartite digraph with parts $X=\bigcup_{t_i \in {\rm in}_T(t)}\{t_i\}\times \Omega$ and $Y=\{t\}\times \Omega$, where the $\delta$-classes in $X$ are the sets $\{t_i\}\times \Omega$, $t_i \in {\rm in}_T(t)$. So the $\delta$-classes  have cardinality $M$, and this is the cardinality of the part $Y$.  Hence the alternets are all isomorphic to a fixed bipartite digraph $\Sigma(m,M)$ and this digraph belongs to  the class $\mathcal{C}$.  In particular, $\Sigma(m,m)$ is the (directed) complete bipartite digraph $\overrightarrow{K}_{m,m}$, and $D(m,m)$ is the digraph with vertex set $\mathbb{Z}\times \mathbb{Z}_m$ and edge set $\{((i,x),(i+1,y))\mid i \in \mathbb{Z}, x,y \in \mathbb{Z}_m\}$.

Note that if $M<\aleph_0$, the digraph $D(m,M)$ is locally finite, and if $M=\aleph_0$, $D(m,M)$ has infinite in-valency.
Moreover, for $x \in D$ and $i \geq 1$, ${\rm desc}^i(x)$ has cardinality $m$ and the digraph induced on ${\rm desc}^i(x) \cup {\rm desc}^{i+1}(x)$ is the (directed) complete bipartite digraph $\overrightarrow{K}_{m,m}$}. 

\end{example}

 Let $m$ be an integer greater than or equal to 1 and let $M$ be either an integer with $M\geq m$, or let $M=\aleph_0$. Throughout  the notation $\Sigma(m,M)$ will be used for the fixed bipartite digraph in $\mathcal{C}$ to which all the alternets of $D(m,M)$ are isomorphic.

\begin{proposition}  \label{M} Let $D$ be an edge transitive, weakly descendant-homogeneous digraph of finite out-valency $m>1$ with property $Z$.  Suppose the alternets of $D$ are isomorphic to some $\Sigma  \in \mathcal{C}$. Then the relation $\delta$ is non-trivial and both ${\rm Al}(D)$ and  the quotient digraph $D/\delta$  are isomorphic to an infinite regular tree $T$ of out-valency 1 and in-valency $k$, where  $k$ is the number of $\delta$-classes in  $\Sigma$. 
\end{proposition}
\begin{proof} Let $R$ be as in Definition \ref{Rclass}. Let $A,B$ be alternets in $D$ with $(A,B)$ an edge in ${\rm Al}(D)$.  Let $U$ denote the $R$-class $Y_A \cap X_B$.
By Lemma \ref{alternets}, $|Y_A| \geq | U| \geq  |{\bf b}|$, where $ {\bf b}$  is a $\delta$-class contained in $U$.  Now assume further that  $A$ is a member of $\mathcal{C}$. Then $Y_A$ and ${\bf b}$ have the same cardinality and, since ${\bf b} \subseteq U \subseteq Y_A $, we conclude that $Y_A=  U=  {\bf b}$. This means that

 (i) $Y_A$ is contained in $X_B$,  and 
 
 (ii)  each  $\delta$-class in $X_B$ is the set of sinks of some alternet.
 
\noindent It follows from (i)  that ${\rm Al}(D)$ has out-valency 1 and is therefore a tree. From (ii) we know that the in-valency  in ${\rm Al}(D)$ is equal to $k$, where  $k$ is the number of $\delta$-classes in  $X_B$. By edge transitivity, we conclude that ${\rm Al}(D)$ is an infinite regular tree of out-valency 1 and in-valency $k$. 

We now show that the quotient digraph $D/\delta$ is also isomorphic to $T$. Let $\bf{u}$ be a $\delta$-class in $D$. Then $\bf{u}$ is contained in the set $X_A$ of sources of some alternet $A$. So ${\rm desc}^1(x) \subseteq Y_A$ for all $x \in \bf{u}$. As shown above, $Y_A=\bf{b}$ for some   $\delta$-class $\bf{b}$ in $D$, and it follows from this fact that for all $x,y \in \bf{u}$ and edges $(x,v),(y,w)$ in $D$,  we have ${\bf{v}=\bf{w}}=Y_A$. This proves that $D/\delta$ has out-valency 1. Finally, since the number of $\delta$-classes in $X_A$ is equal to $k$, we conclude that  $D/\delta$ has in-valency $k$. 
 \end{proof}
 
 \vspace{0.2cm}
 
 The next result  gives a general description of the distance-transitive digraphs $D$ for which the reachability relation is not universal and the alternets belong to the class $\mathcal{C}$.


\begin{theorem}\label{description} Let $D$ be a distance-transitive, weakly descendant-homogeneous digraph of finite out-valency $m>1$ with no directed cycles. Suppose the reachability relation is not universal and the alternets of $D$ belong to $ \mathcal{C}$. Then the relation   $\delta$ is non-trivial and the quotient digraph $D/\delta$ is a tree of out-valency 1 and in-valency $k$, where  $k$ is the number of $\delta$-classes in  an alternet. 
In particular, if  the alternets in $D$ are isomorphic to $\Sigma(m,M)$, for some integer $M \geq m$ or $M=\aleph_0$,  then $D$ is isomorphic to  the digraph $D(m,M)$.
 \end{theorem}
 \begin{proof}
By Theorem  \ref{equivalence1} we know  that $D$ has property $Z$, since we are assuming that the reachability relation is not universal. Then, by Proposition \ref{M}, the quotient digraph $D/\delta$ is a tree of out-valency 1 and in-valency $k$, where  $k$ is the number of $\delta$-classes in  $\Sigma$.  In particular, if  there exists  $M$  such that the alternets in $D$ are isomorphic to $\Sigma(m,M)$ then $D$ is isomorphic to $D(m,M)$.
 \end{proof}

 \begin{corollary} Let $D$ be an edge transitive, weakly descendant-homogeneous digraph of finite out-valency $m>1$ with property $Z$.  Suppose that the relation $\delta$ is non-trivial and that the set of sinks of an alternet $A$ has cardinality equal to a prime $p$. Then $A \in \mathcal{C}$ and both ${\rm Al}(D)$ and  the quotient digraph $D/\delta$  are isomorphic to an infinite regular tree $T$ of out-valency 1 and in-valency k, where  $k$ is the number of $\delta$-classes in  $A$. 
 In particular, the in-valency $r$ of $D$ is a multiple of $m$. 
 \end{corollary}
 \begin{proof} Let $A$ be an alternet in $D$. By Proposition  \ref{M},  it is enough to show that $A \in \mathcal{C}$. Let $X$ and $Y$ be the set of sources and the set of sinks of $A$, respectively. Let $R$ be as in Definition \ref{Rclass}, and let $U$ be a $R$-class in $Y$. We know that the  cardinality  of $U$  divides the cardinality of $Y$ and, therefore, is equal to 1 or $p$. By Lemma \ref{alternets}, $U$  is a union of $\delta$-classes and, by assumption,  such classes contain at least two vertices. So the cardinality of $U$ is  $p$ and, therefore,  $U$ consists of a unique $\delta$-class. Hence $U=Y$ and the $\delta$-classes have $p$ vertices. Thus  $A \in \mathcal{C}$. 
  \end{proof}

 \begin{corollary}  \label{M2} Let $D$ be an edge transitive, weakly descendant-homogeneous digraph of finite out-valency $m>1$ with property $Z$.  Suppose  the alternets of $D$ are isomorphic to $\Sigma(m,M)$, for some integer $M \geq m$ or $M=\aleph_0$. Then $D$ is isomorphic to $D(m,M)$. 
In particular, $D$ is highly-arc-transitive. 
\end{corollary}
\begin{proof} It is enough to note that the relation $\delta$ is non-trivial in $D(m,M)$ and, by the construction, the quotient $D(m,M)/\delta$ is an infinite regular tree of out-valency 1 and in-valency equal to the number of $\delta$-classes in $\Sigma(m,M)$. 
  \end{proof}
  
  \vspace{0.3cm}
  
 More generally, we have the following.
 
 \begin{theorem} \label{HAT}Let $D$ be an edge transitive, weakly descendant-homogeneous digraph of finite out-valency $m>1$ with property $Z$.   Suppose the alternets of $D$ are isomorphic to some $\Sigma  \in \mathcal{C}$. Then $D$ is highly-arc-transitive.
\end{theorem}
\begin{proof} We use induction on $s$ to  prove that $D$ is $s$-arc transitive  for all $s \geq 0$. For  $s\leq 1$ the result holds by edge transitivity.  Let $s\geq 2$ and assume $D$ is $s$-arc transitive. To prove that $D$ is $(s+1)$-arc transitive, we show that for any $s$-arc $\gamma$ in $D$, the stabiliser $G_{\gamma}$ is transitive on the set ${\rm desc}^1(u)$ of out-neighbours of the end vertex $u$ of $\gamma$. 

By Proposition  \ref{M},  the relation $\delta$ is non-trivial and the quotient digraph $D/\delta$ is  isomorphic to an infinite regular tree $T$ of out-valency 1 and in-valency $k$. Consider the 
 $\delta$-class $\bf{u}$ containing $u$ and let $U$ denote the finitely generated subdigraph induced on ${\rm desc}(\bf{u})$. It follows from  edge transitivity of $D$ that for any $a,b$ in ${\rm desc}^1(u)$ there is an automorphism $h$ of ${\rm desc}(u)$ which takes $a$ to $b$.  Since the subdigraph on ${\bf{u}} \cup {\rm desc}^1(\bf{u})$ is complete bipartite, $h$ extends to an automorphism $h'$ of $U$ which fixes $\bf{u}$ pointwise. By weak descendant-homogeneity, and the fact that $D/\delta$ is a tree, $h'$ extends to an automorphism $g$ of $D$ such that 
 
 (i) $g$ agrees with $h'$ on $\bf{u} \cup {\rm desc}(\bf{u})$, and 
 
 (ii) $g$ fixes pointwise the set ${\rm anc}({\bf{u}})$ of ancestors of {\bf{u}}, that is, the set of vertices $ x \in D$ for which there is an $s$-arc from $x$ to a vertex $y \in \bf{u}$. \\
 Then  from (i) we  know that $g$ takes $a$ to $b$. Since $\gamma$ is contained in ${\rm anc}({\bf{u}})$, we conclude from (ii)  that $g \in G_{\gamma}$. The result follows. 
  \end{proof} 
  
  \vspace{0.3cm}

 In particular, we now consider the case when the digraph ${\rm Al}(D)$ is isomorphic to $Z$.  
 
 \begin{lemma}  \label{alternetZ} Let $D$ be an edge transitive  digraph of finite out-valency $m>1$ such that ${\rm Al}(D)$  is isomorphic to $Z$. Then $D$ is two-ended locally finite digraph with equal in- and out-valency and with  property $Z$. 
\end{lemma}
\begin{proof} 
  Let $F$ be a finite set of vertices of $D$. We show that the digraph obtained from $D$ by removing $F$ contains at two infinite components. Since ${\rm Al}(D)$  is isomorphic to $Z$ we can enumerate the alternets of $D$ as $\{A_i, i \in \mathbb{Z}\}$ with $(A_i,A_{i+1})$ an edge for all $i$. As $F$ is finite, $F$ is a subset of the vertex set of a union of finitely many alternets. By vertex transitivity of $D$, we may assume $F \subseteq \bigcup_{i=1}^n VA_i$, where $VA_i$ denotes the set of vertices of the alternet $A_i$.
 The digraph obtained from ${\rm Al}(D)$ by deleting the vertices in $A_1, \ldots, A_n$ contains two infinite components, namely the two infinite arcs $\ldots, A_{-i}, \ldots, A_0$
 and $A_{n+1},A_{n+2} \ldots$. Since each alternet $A_i$ is connected, both the subdigraphs of $D$ on $\bigcup_{j \leq 0} VA_j$ and on $\bigcup_{j\geq n+1} VA_i$ are connected. Hence $D \setminus F$ contains two infinite components.
   
Let $X$ be the underlying graph of $D$. Then $X$   is a vertex transitive  two-ended graph. It follows from \cite[Theorem 7]{diestel} that  $X$ is locally finite. Therefore $D$ is locally finite. The fact that $D$ has equal in- and out-valency now follows from \cite[Corollary 9(c)]{Moller5}. 
Finally,  $D$ has property $Z$ since ${\rm Al}(D)$ has property $Z$.
  \end{proof} 
  
 \vspace{0.3cm} 

\begin{proposition} \label{hatZ}  Let $D$ be a distance-transitive, weakly descendant-homogeneous digraph of finite out-valency $m>1$ with no directed cycles. If  ${\rm Al}(D)$ has in-valency 1 then $D$ is a two-ended locally finite digraph with equal in- and out-valency and with  property $Z$. \end{proposition}
\begin{proof}  
By Lemma \ref{alternetZ}, it is enough to show that ${\rm Al}(D)$ has out-valency 1. Seeking a contradiction, assume ${\rm Al}(D)$ has out-valency greater than 1. Let $A,B \in {\rm Al}(D)$ such that $(A,B)$ is an edge.  Since $A$ is connected, there is a vertex $x \in X_A$ and an alternet $C \neq B$ such that  $(A,C)$ is an edge with  $\Gamma(x) \cap X_B \neq \varnothing$ and $\Gamma(x) \cap X_C \neq \varnothing$. Now,  as ${\rm Al}(D)$ is a tree, the descendant sets of $B$ and $C$ do not intersect. So  for any vertex $b$ in $X_B$ (the set of sources of the alternet $B$) and any vertex $c$ in $X_C$ (the sources of the alternet $C$), the intersection ${\rm desc}(b)\cap {\rm desc}(c)$ in $D$ is empty. In particular,  ${\rm desc}(u)\cap {\rm desc}(v)=\varnothing$ for  $u\in \Gamma(x) \cap X_B$ and $v\in \Gamma(x) \cap X_C$. It then follows from Lemma \ref{propZ1}  that $D$ does not have property Z. \end{proof}

 \begin{corollary}   Let $D$ be a highly-arc-transitive, weakly descendant-homogeneous digraph of finite out-valency $m>1$. If  ${\rm Al}(D)$ has in-valency 1 then $D$ is the digraph induced on the descendant set of a directed line in $D$.
\end{corollary}
 \begin{proof} By Proposition  \ref{hatZ}, $D$ is a two-ended digraph. The result then follows from  Theorem  \ref{2end}.  \end{proof}

\section{Final Remarks}

Let $D$ be a  distance-transitive, weakly descendant-homogeneous digraph of finite out-valency $m>1$  with property $Z$. By Theorem \ref{primitive},  there is a natural number $n\geq 1$ such that  the relation $\delta_n$ is non-trivial in $D$.  The digraphs in Example \ref{example1} are digraphs with these properties.  
Note that the relation $\delta_1$ is non-trivial in the digraph  $D(m,M)$ with $(m,M)\neq (1,1)$, and  $\delta_1$ is trivial in the digraph $D(1,1)$ since this is isomorphic to $Z$.  There are examples of distance-transitive digraphs with property $Z$ in which   $\delta_1$  is trivial and $\delta_n$ is non-trivial for some $n \geq 2$ (see Section 4.2 of  {\cite{Amato2}}), but they are not weakly descendant-homogeneous. This suggests the following question.

 \begin{question} Let $D$ be a distance-transitive, weakly descendant-homogeneous digraph of finite out-valency $m>1$ with property $Z$. Is the relation $\delta_1$  non-trivial in $D$? \end{question}
 
 Also,  the following question arises.
 
  \begin{question}\label{Q2} Let $D$ be a distance-transitive, weakly descendant-homogeneous digraph of finite out-valency $m>1$ with property $Z$. Suppose $D$ is not a tree. If the alternets of $D$ are finite, do they belong to  the class $ \mathcal{C}$?
 \end{question}

If the answer to Question \ref{Q2} is positive, then, together with Theorem  \ref{HAT}, it would prove the following conjecture.

 \begin{conjecture} Let $D$ be a distance-transitive, weakly descendant-homogeneous digraph of finite out-valency $m>1$ with property $Z$. If the alternets of $D$ are finite then $D$ is highly-arc-transitive. 
 \end{conjecture}

 \section*{Acknowledgement}
 I would like to thank the anonymous referee for helpful comments and suggestions which improved the quality of the paper. I would also like to thank Professor David Evans for the critical reading of the manuscript.

 

\end{document}